\documentclass{amsproc}
\usepackage{amsopn}
\usepackage{amssymb, amscd}

\newtheorem{theorem}{Theorem}[section]
\newtheorem{lemma}[theorem]{Lemma}

\theoremstyle{definition}
\newtheorem{definition}[theorem]{Definition}
\newtheorem{example}[theorem]{Example}

\theoremstyle{remark}
\newtheorem{remark}[theorem]{Remark}

\numberwithin{equation}{section}

\newcommand{\nc}{\newcommand}

\nc{\fg}{\mathfrak{f} } \nc{\vg}{\mathfrak{v} } \nc{\wg}{\mathfrak{w} }
\nc{\zg}{\mathfrak{z} } \nc{\ngo}{\mathfrak{n} } \nc{\kg}{\mathfrak{k} }
\nc{\mg}{\mathfrak{m} } \nc{\bg}{\mathfrak{b} } \nc{\ggo}{\mathfrak{g} }
\nc{\ggob}{\overline{\mathfrak{g}} } \nc{\sog}{\mathfrak{so} }
\nc{\sug}{\mathfrak{su} } \nc{\spg}{\mathfrak{sp} } \nc{\slg}{\mathfrak{sl} }
\nc{\glg}{\mathfrak{gl} } \nc{\cg}{\mathfrak{c} } \nc{\rg}{\mathfrak{r} }
\nc{\hg}{\mathfrak{h} } \nc{\tg}{\mathfrak{t} } \nc{\ug}{\mathfrak{u} }
\nc{\dg}{\mathfrak{d} } \nc{\ag}{\mathfrak{a} } \nc{\pg}{\mathfrak{p} }
\nc{\sg}{\mathfrak{s} } \nc{\pca}{\mathcal{P}} \nc{\nca}{\mathcal{N}}
\nc{\lca}{\mathcal{L}} \nc{\oca}{\mathcal{O}} \nc{\mca}{\mathcal{M}}
\nc{\tca}{\mathcal{T}} \nc{\aca}{\mathcal{A}} \nc{\cca}{\mathcal{C}}
\nc{\gca}{\mathcal{G}} \nc{\sca}{\mathcal{S}} \nc{\hca}{\mathcal{H}}
\nc{\bca}{\mathcal{B}} \nc{\dca}{\mathcal{D}} \nc{\val}{\operatorname{val}}

\nc{\vp}{\varphi} \nc{\ddt}{\tfrac{{\rm d}}{{\rm d}t}} \nc{\im}{\mathtt{i}}

\nc{\SO}{\mathrm{SO}} \nc{\Spe}{\mathrm{Sp}} \nc{\Sl}{\mathrm{SL}}
\nc{\SU}{\mathrm{SU}} \nc{\Or}{\mathrm{O}} \nc{\U}{\mathrm{U}} \nc{\Gl}{\mathrm{GL}}
\nc{\Se}{\mathrm{S}} \nc{\Cl}{\mathrm{Cl}} \nc{\Spein}{\mathrm{Spin}}
\nc{\Pin}{\mathrm{Pin}} \nc{\G}{\mathrm{GL}_n(\RR)} \nc{\g}{\mathfrak{gl}_n(\RR)}

\nc{\RR}{{\Bbb R}} \nc{\HH}{{\Bbb H}} \nc{\CC}{{\Bbb C}} \nc{\ZZ}{{\Bbb Z}}
\nc{\FF}{{\Bbb F}} \nc{\NN}{{\Bbb N}} \nc{\QQ}{{\Bbb Q}} \nc{\PP}{{\Bbb P}}

\nc{\vs}{\vspace{.2cm}} \nc{\vsp}{\vspace{1cm}} \nc{\ip}{\langle\cdot,\cdot\rangle}
\nc{\ipp}{(\cdot,\cdot)} \nc{\la}{\langle} \nc{\ra}{\rangle} \nc{\unm}{\tfrac{1}{2}}
\nc{\unc}{\tfrac{1}{4}} \nc{\und}{\tfrac{1}{16}} \nc{\no}{\vs\noindent}
\nc{\lam}{\Lambda^2(\RR^n)^*\otimes\RR^n} \nc{\tangz}{{\rm T}^{\rm Zar}}
\nc{\nor}{{\sf n}}  \nc{\mum}{/\!\!/} \nc{\kir}{/\!\!/\!\!/}
\nc{\Ri}{\tfrac{4\Ric_{\mu}}{||\mu||^2}} \nc{\ds}{\displaystyle}
\nc{\ben}{\begin{enumerate}} \nc{\een}{\end{enumerate}} \nc{\f}{\frac}
\nc{\lb}{[\cdot,\cdot]} \nc{\isn}{\tfrac{1}{||v||^2}}

\nc{\He}{\operatorname{Hess}} \nc{\ad}{\operatorname{ad}}
\nc{\Ad}{\operatorname{Ad}} \nc{\rank}{\operatorname{rank}}
\nc{\Irr}{\operatorname{Irr}} \nc{\End}{\operatorname{End}}
\nc{\Aut}{\operatorname{Aut}} \nc{\Inn}{\operatorname{Inn}}
\nc{\Der}{\operatorname{Der}} \nc{\Ker}{\operatorname{Ker}}
\nc{\Iso}{\operatorname{I}} \nc{\Diff}{\operatorname{D}} \nc{\Lie}{\operatorname{L}}
\nc{\tr}{\operatorname{tr}} \nc{\dif}{\operatorname{d}}
\nc{\sen}{\operatorname{sen}} \nc{\modu}{\operatorname{mod}}
\nc{\Ric}{\operatorname{R}} \nc{\Ricci}{\operatorname{Ric}}
\nc{\sym}{\operatorname{sym}} \nc{\symac}{\operatorname{sym^{ac}}}
\nc{\symc}{\operatorname{sym^{c}}} \nc{\scalar}{\operatorname{sc}}
\nc{\grad}{\operatorname{grad}} \nc{\ricci}{\operatorname{ric}}
\nc{\ricciac}{\operatorname{ric^{ac}}} \nc{\riccic}{\operatorname{ric^{c}}}
\nc{\riccig}{\operatorname{ric^{\gamma}}} \nc{\Rin}{\operatorname{M}}
\nc{\Le}{\operatorname{L}} \nc{\tang}{\operatorname{T}}
\nc{\level}{\operatorname{level}} \nc{\rad}{\operatorname{r}}
\nc{\abel}{\operatorname{ab}} \nc{\CH}{\operatorname{CH}}
\nc{\mcc}{\operatorname{mcc}} \nc{\Adj}{\operatorname{Adj}}

\begin{document}

\title{Einstein solvmanifolds and nilsolitons}

\author{Jorge Lauret}

\address{FaMAF and CIEM, Universidad Nacional de C\'ordoba, C\'ordoba, Argentina}
\email{lauret@mate.uncor.edu}
\thanks{Supported by grants from CONICET, SECYT (Univ. Nac. C\'ordoba) and FONCYT}

\date{May 15, 2008}

\dedicatory{Dedicated to Isabel Dotti and Roberto Miatello on the occasion of their
60th birthday.}


\begin{abstract}
The purpose of the present expository paper is to give an account of the recent
progress and present status of the classification of solvable Lie groups admitting
an Einstein left invariant Riemannian metric, the only known examples so far of
noncompact Einstein homogeneous manifolds.  The problem turns to be equivalent to
the classification of Ricci soliton left invariant metrics on nilpotent Lie groups.
\end{abstract}

\maketitle

\tableofcontents

\section{Introduction}\label{intro}

Let $M$ be a differentiable manifold.  The question of whether there is a `best'
Riemannian metric on $M$ is intriguing.  A great deal of deep results in
Riemannian geometry have been motivated, and even inspired, by this single natural
question. For several good reasons, an Einstein metric is a good candidate, if
not the best, at least a very distinguished one (see \cite[Chapter 0]{Bss}).  A
Riemannian metric $g$ on $M$ is called {\it Einstein} if its Ricci tensor $\ricci_g$
satisfies
\begin{equation}\label{eeq}
\ricci_g=cg, \qquad \mbox{for some}\; c\in\RR.
\end{equation}

This notion can be traced back to \cite{Hlb}, where Einstein metrics emerged as
critical points of the total scalar curvature functional on the space of all metrics
on $M$ of a given volume. Equation (\ref{eeq}) is a non-linear second order PDE
(recall that the number of parameters is $\tfrac{n(n+1)}{2}$ on both sides,
$n=\dim{M}$), which also gives rise to some hope, but a good understanding of the
solutions in the general case seems far from being attained.  A classical reference for Einstein
manifolds is the book \cite{Bss}, and some updated expository articles are
\cite{And}, \cite{LbrWng}, \cite[III,C.]{Brg1} and \cite[11.4]{Brg2}.

The Einstein condition (\ref{eeq}) is very subtle, even when restricted to almost
any subclass of metrics on $M$ one may like.  It is too strong to allow general
existence results, and sometimes even just to find a single example, and at the same
time, it is too weak to get obstructions or classification results.

But maybe the difficulty comes from PDEs, so let us `algebrize' the problem
(algebra is always easier for a geometer ...).  Let us consider homogeneous
Riemannian manifolds.  Indeed, the Einstein equation for a homogeneous metric is
just a system of $\tfrac{n(n+1)}{2}$ algebraic equations, but unfortunately, a quite
involved one, and the following main general question is still open in both compact
and noncompact cases:

\begin{quote}
Which homogeneous spaces $G/K$ admit a $G$-invariant Einstein Riemannian metric?
\end{quote}

We refer to \cite{BhmWngZll} and the references therein for an update in the compact
case.  In the noncompact case, the only known examples until now are all of a very
particular kind; namely, simply connected solvable Lie groups endowed with a left
invariant metric (so called {\it solvmanifolds}). According to the following long
standing conjecture, these might exhaust all the possibilities for noncompact
homogeneous Einstein manifolds.

\begin{quote}
{\bf Alekseevskii's conjecture} \cite[7.57]{Bss}.  If $G/K$ is a homogeneous
Einstein manifold of negative scalar curvature then $K$ is a maximal compact
subgroup of $G$ (which implies that $G/K$ is a solvmanifold when $G$ is a linear
group).
\end{quote}

The conjecture is wide open, and it is known to be true only for $\dim \leq 5$, a
result which follows from the complete classification in these dimensions given in
\cite{Nkn}.  One of the most intriguing facts related to this conjecture, and maybe
the only reason so far to consider Alekseevskii's conjecture as too optimistic, is
that the Lie groups $\Sl_n(\RR)$, $n\geq 3$, do admit left invariant metrics of
negative Ricci curvature, as well as does any complex simple Lie group (see
\cite{DttLt}, \cite{DttLtMtl}).  However, an inspection of the eigenvalues of the
Ricci tensors in these examples shows that they are far from being close to each
other, giving back some hope.

Let us now consider the case of left invariant metrics on Lie groups.  Let $\ggo$ be
a real Lie algebra.  Each basis $\{ X_1,...,X_n\}$ of $\ggo$ determines structural
constants $\{ c_{ij}^k\}\subset\RR$ given by
$$
[X_i,X_j]=\sum_{k=1}^nc_{ij}^kX_k, \qquad 1\leq i,j\leq n.
$$
The left invariant metric on any Lie group with Lie algebra $\ggo$ defined by the
inner product given by $\la X_i,X_j\ra=\delta_{ij}$ is Einstein if and only if the
$\tfrac{n^2(n+1)}{2}$ numbers $c_{ij}^k$'s satisfy the following $\tfrac{n(n+1)}{2}$
algebraic equations for some $c\in\RR$:
\begin{equation}\label{condi}
\sum_{kl}-\unm c_{ik}^lc_{jk}^l +\unc c_{kl}^ic_{kl}^j-\unm c_{jk}^lc_{il}^k+\unm
c_{kl}^lc_{ki}^j +\unm c_{kl}^lc_{kj}^i = c\delta_{ij}, \quad 1\leq i\leq j\leq n.
\end{equation}

In view of this, one may naively think that the classification of Einstein left
invariant metrics on Lie groups is at hand.  However, the following natural
questions remain open:

\begin{itemize}
\item[ ]
\item[(i)] Is any Lie group admitting an Einstein left invariant metric either solvable or compact?

\item[ ]\item[(ii)] Does every compact Lie group admit only finitely many Einstein left invariant metrics up to isometry and scaling?

\item[ ]\item[(iii)]  Which solvable Lie groups admit an Einstein left invariant metric?

\item[ ]
\end{itemize}

We note that question (i) is just Alekseevskii Conjecture restricted to Lie groups,
and question (ii) is contained in \cite[7.55]{Bss}.  The only group for which the
answer to (ii) is known is $\SU(2)$, where there is only one (see \cite{Mln}).  For
most of the other compact simple Lie groups many Einstein left invariant metrics
other than minus the Killing form are explicitly known (see \cite{DtrZll}).

Even if one is very optimistic and believes that Alekseevskii Conjecture is true, a
classification of Einstein metrics in the noncompact homogeneous case will depend on
some kind of answer to question (iii).  The aim of this expository paper is indeed
to give a report on the present status of the study of Einstein solvmanifolds.

Perhaps the main difficulty in trying to decide if a given Lie algebra $\ggo$ admits
an Einstein inner product is that one must check condition (\ref{condi}) for any
basis of $\ggo$, and there are really too many of them.  In other words, there are
too many left invariant metrics on a given Lie group, any inner product on the
vector space $\ggo$ is playing. This is quite in contrast to what happens in
homogeneous spaces $G/K$ with not many different $\Ad(K)$-irreducible components in
the decomposition of the tangent space $\tang_{eK}(G/K)$.  Another obstacle is how to
recognize your Lie algebra by just looking at the structural constants $c_{ij}^k$'s.
Though even in the case when we have two solutions to (\ref{condi}), and we know they define
the same Lie algebra, to be able to guarantee that they are not isometric, i.e. that we
really have two Einstein metrics, is usually involved.

If we fix a basis $\{ X_1,...,X_n\}$ of $\ggo$, then instead of varying all possible
sets of structural constants $\{ c_{ij}^k\}$'s by running over all bases, one may
act on the Lie bracket $\lb$ by $g.\lb=g[g^{-1}\cdot,g^{-1}\cdot]$, for any
$g\in\Gl(\ggo)$, and look at the structural constants of $g.\lb$ with respect to the
fixed basis $\{ X_1,...,X_n\}$.  This give rises to an orbit $\Gl(\ggo).\lb$ in the
vector space $V:=\Lambda^2\ggo^*\otimes\ggo$ of all skew-symmetric bilinear maps
from $\ggo\times\ggo$ to $\ggo$, which parameterizes, from a different point of
view, the set of all inner products on $\ggo$. Indeed, if $\ip$ is the inner product
defined by $\la X_i,X_j\ra=\delta_{ij}$ then

\begin{quote}
$(\ggo,g.\lb,\ip)$ is isometric to $(\ggo,\lb,\la g\cdot,g\cdot\ra)$ for any
$g\in\Gl(\ggo)$.
\end{quote}

The subset $\lca\subset V$ of those elements satisfying the Jacobi condition is
algebraic, $\Gl(\ggo)$-invariant and the $\Gl(\ggo)$-orbits in $\lca$ are precisely
the isomorphism classes of Lie algebras.  $\lca$ is called the {\it variety of Lie
algebras}. Furthermore, if $\Or(\ggo)\subset\Gl(\ggo)$ denotes the subgroup of
$\ip$-orthogonal maps, then two points in $\Gl(\ggo).\lb$ which lie in the same
$\Or(\ggo)$-orbit determine isometric left invariant metrics, and the converse holds
if $\ggo$ is completely solvable (see \cite{Alk2}).

This point of view is certainly a rather tempting invitation to try to use geometric
invariant theory in any problem which needs a running over all left invariant
metrics on a given Lie group, or even on all Lie groups of a given dimension.  We
shall see throughout this article that indeed, starting in \cite{Hbr}, the approach `by varying
Lie brackets' has been very fruitful in the study of Einstein solvmanifolds
during the last decade.

The latest fashion generalization of Einstein metrics, although they were introduced
by R. Hamilton more than twenty years ago, is the notion of {\it Ricci soliton}:
\begin{equation}\label{rseq}
\ricci_g=cg+L_Xg, \qquad\mbox{for some}\; c\in\RR, \quad X\in\chi(M),
\end{equation}

\noindent where $L_Xg$ is the usual Lie derivative of $g$ in the direction of the
field $X$.  A more intuitive equivalent condition to (\ref{rseq}) is that $\ricci_g$
is tangent at $g$ to the space of all metrics which are homothetic to $g$ (i.e.
isometric up to a constant scalar multiple).  Recall that Einstein means
$\ricci_{g}$ tangent to $\RR_{>0}g$.  Ricci solitons correspond to solutions of the
Ricci flow
$$
\ddt g(t)=-2\ricci_{g(t)},
$$
that evolves self similarly, that is, only by scaling and the action by
diffeomorphisms, and often arise as limits of dilations of singularities of the
Ricci flow.  We refer to \cite{soliton}, \cite{GntIsnKnp}, \cite{libro} and the
references therein for further information on the Hamilton-Perelman theory of Ricci
flow and Ricci solitons and the role played by nilpotent Lie groups in the story.

A remarkable fact is that if $S$ is an Einstein solvmanifold, then the metric
restricted to the submanifold $N:=[S,S]$ is a Ricci soliton, and conversely, any
Ricci soliton left invariant metric on a nilpotent Lie group $N$ (called {\it
nilsolitons}) can be uniquely `extended' to an Einstein solvmanifold.  This
one-to-one correspondence is complemented with the uniqueness up to isometry of
nilsolitons, which finally turns the classification of Einstein solvmanifolds into a
classification problem on nilpotent Lie algebras.  These are not precisely good
news.  Historically, as the literature and experience shows us, any classification
problem involving nilpotent Lie algebras is simply a headache.

\vs \noindent {\it Acknowledgements.}  I am grateful to the Scientific
Committee for the invitation to give a talk at the `Sixth Workshop on Lie Theory and
Geometry', November 13-17, 2007, Cruz Chica, C\'ordoba, Argentina.  I wish to thank Yuri Nikolayevsky and Cynthia Will for very useful comments on a first version of the paper, and to Roberto Miatello for going over the manuscript.  I also wish to express my gratitude to the young collaborators Alejandra \'Alvarez, Adri\'an
Andrada, Gast\'on Garc\'{\i}a and Emilio Lauret for the invaluable help they generously provided to the organization of the workshop.

\section{Structure and uniqueness results on Einstein solvmanifolds}\label{pre}

A {\it solvmanifold} is a simply connected solvable Lie group $S$ endowed with a
left invariant Riemannian metric.  A left invariant metric on a Lie group $G$ will
be always identified with the inner product $\ip$ determined on the Lie algebra
$\ggo$ of $G$, and the pair $(\ggo,\ip)$ will be referred to as a {\it metric Lie
algebra}.  If $S$ is a solvmanifold and $(\sg,\ip)$ is its metric solvable Lie
algebra, then we consider the $\ip$-orthogonal decomposition
$$
\sg=\ag\oplus\ngo,
$$
where $\ngo:=[\sg,\sg]$ is the derived algebra (recall that $\ngo$ is nilpotent).

\begin{definition}\label{stan}
A solvmanifold $S$ is said to be {\it standard} if
$$
[\ag,\ag]=0.
$$
\end{definition}

This is a very simple algebraic condition, which may appear as kind of technical,
but it has nevertheless played an important role in many questions in homogeneous
Riemannian geometry:

\begin{itemize}
\item \cite{GndPttVnb} K$\ddot{{\rm a}}$hler-Einstein noncompact homogeneous manifolds are all standard solvmanifolds.

\item\cite{Alk,Crt} Every quaternionic K$\ddot{{\rm a}}$hler solvmanifold (completely real) is standard.

\item \cite{AznWls} Any homogeneous manifold of nonpositive sectional curvature is a standard solvmanifold.

\item \cite{Hbr2} All harmonic noncompact homogeneous manifolds are standard solvmanifolds (with $\dim{\ag}=1$).
\end{itemize}

Partial results on the question of whether Einstein solvmanifolds are all standard
were obtained for instance in \cite{Hbr} and \cite{Sch}, who gave several sufficient
conditions. The answer was known to be yes in dimension $\leq 6$ (see \cite{NktNkn})
and followed from a complete classification of Einstein solvmanifolds in these
dimensions.  On the other hand, it is proved in \cite{Nkl0} that many classes of
nilpotent Lie algebras can not be the nilradical of a non-standard Einstein
solvmanifold.

\begin{theorem}\cite{standard}\label{stand}
Any Einstein solvmanifold is standard.
\end{theorem}

An idea of the proof of this theorem will be given in Section \ref{proof}.  Standard
Einstein solvmanifolds were extensively investigated in \cite{Hbr}, where the
remarkable structural and uniqueness results we next describe are derived.  Recall
that combined with Theorem \ref{stand}, all of these results are now valid for any
Einstein solvmanifold.

\begin{theorem}\label{u}\cite[Section 5]{Hbr} {\rm (}{\bf Uniqueness}{\rm )}
A simply connected solvable Lie group admits at most one standard Einstein left
invariant metric up to isometry and scaling.
\end{theorem}

A more general result is actually valid: if a noncompact homogeneous manifold $G/K$
with $K$ maximal compact in $G$ admits a $G$-invariant metric $g$ isometric to an
Einstein solvmanifold, then $g$ is the unique $G$-invariant Einstein metric on $G/K$
up to isometry and scaling.  This is in contrast to the compact homogeneous case,
where many pairwise non isometric $G$-invariant Einstein metrics might exist (see
\cite{BhmWngZll} and the references therein), although it is open if only finitely
many (see \cite[7.55]{Bss}).

In the study of Einstein homogeneous manifolds, the compact case is characterized by
the positivity of the scalar curvature and Ricci flat implies flat (see
\cite{AlkKml}).  The following conditions on an Einstein solvmanifold $S$ are
equivalent:

\begin{itemize}
\item[(i)] $\sg$ is unimodular (i.e. $\tr{\ad{X}}=0$ for all $X\in\sg$).

\item[(ii)] $S$ is Ricci flat (i.e. $\scalar(S)=0$).

\item[(iii)] $S$ is flat.
\end{itemize}

We can therefore consider from now on only nonunimodular solvable Lie algebras.

\begin{theorem}\label{ror}\cite[Section 4]{Hbr} {\rm (}{\bf Rank-one reduction}{\rm )}
Let $\sg=\ag\oplus\ngo$ be a nonunimodular solvable Lie algebra endowed with a
standard Einstein inner product $\ip$, say with $\ricci_{\ip}=c\ip$.  Then $c<0$
and, up to isometry, it can be assumed that $\ad{A}$ is symmetric for any $A\in\ag$.
In that case, the following conditions hold.
\begin{itemize}
\item[(i)] There exists $H\in\ag$ such that the eigenvalues of $\ad{H}|_{\ngo}$ are all positive integers without a common divisor.

\item[(ii)] The restriction of $\ip$ to the solvable Lie algebra $\RR H\oplus\ngo$ is also Einstein.

\item[(iii)] $\ag$ is an abelian algebra of symmetric derivations of $\ngo$ and the inner product on $\ag$ must be given by $\la A,A'\ra=-\tfrac{1}{c}\tr{\ad{A}\ad{A'}}$ for all $A,A'\in\ag$.
\end{itemize}
\end{theorem}

The Ricci tensor for these solvmanifolds has the following simple formula.

\begin{lemma}\label{fr}
Let $S$ be a standard solvmanifold such that $\ad{A}$ is symmetric and nonzero for
any $A\in\ag$. Then the Ricci tensor of $S$ is given by
\begin{itemize}
\item[(i)] $\ricci(A,A')=-\tr{\ad{A}\ad{A'}}$ for all $A,A'\in\ag$.

\item[(ii)] $\ricci(\ag,\ngo)=0$.

\item[(iii)] $\ricci(X,Y)=\ricci_{\ngo}(X,Y) -\la\ad{H}(X),Y\ra$, for all $X,Y\in\ngo$, where $\ricci_{\ngo}$ is the Ricci tensor of $(\ngo,\ip|_{\ngo\times\ngo})$ and $H\in\ag$ is defined by $\la H,A\ra=\tr{\ad{A}}$ for any $A\in\ag$.
\end{itemize}
\end{lemma}

The natural numbers which have appeared as the eigenvalues of $\ad{H}$ when
$(\sg,\ip)$ is Einstein play a very important role.

\begin{definition}\label{et}
If $d_1,...,d_r$ denote the corresponding multiplicities of the positive integers
without a common divisor $k_1<...<k_r$ given by Theorem \ref{ror}, (i), then the
tuple
$$
(k;d)=(k_1<...<k_r;d_1,...,d_r)
$$
is called the {\it eigenvalue type} of the Einstein solvmanifold $(\sg,\ip)$.
\end{definition}

We find here the first obstruction: if a solvable Lie algebra admits an Einstein
inner product then the nilpotent Lie algebra $\ngo=[\sg,\sg]$ is $\NN$-{\it graded},
that is, there is a decomposition $\ngo=\ngo_1\oplus...\oplus\ngo_r$ such that
$[\ngo_i,\ngo_j]\subset\ngo_{i+j}$ for all $i,j$ (recall that some of the $\ngo_i$'s
might be trivial).  This is precisely the decomposition into eigenspaces of the derivation
with positive integer eigenvalues $\ad{H}$.  Another important consequence of
Theorem \ref{ror} is that to study Einstein solvmanifolds, it will be enough to consider
{\it rank-one} (i.e. $\dim{\ag}=1$) metric solvable Lie algebras, since every higher
rank Einstein solvmanifold will correspond to a unique rank-one Einstein
solvmanifold and certain abelian subalgebra $\ag$ of derivations of $\ngo$
containing $\ad{H}$.  Recall that how to extend the inner product is determined by
Theorem \ref{ror}, (iii).

Let $\mca$ be the moduli space of all the isometry classes of Einstein solvmanifolds
of a given dimension with scalar curvature equal to $-1$, endowed with the
$C^{\infty}$-topology.  Notice that any $n$-dimensional solvmanifold $S$ is
diffeomorphic to the euclidean space $\RR^n$, and so any $S$ can be viewed as a
Riemannian metric on $\RR^n$ (which is in addition invariant by some transitive
solvable Lie group of diffeomorphisms of $\RR^n$).

\begin{theorem}\label{m}\cite[Section 6]{Hbr} {\rm (}{\bf Moduli space}{\rm )}
In every dimension, only finitely many eigenvalue types occur, and each eigenvalue
type $(k;d)$ determines a compact path connected component $\mca_{(k;d)}$ of $\mca$,
homeomorphic to a real semialgebraic set.
\end{theorem}

Results on the topology and `dimension' of the moduli spaces $\mca_{(k;d)}$ near a
rank-one symmetric space are obtained in \cite[Section 6.5]{Hbr}.  This has also
been done for many other symmetric spaces (rank $\geq 2$) in \cite{GrdKrr}, where
even explicit examples are exhibited to describe a neighborhood.  The moduli spaces
$\mca_{(1<2;q,p)}$ are studied in detail in \cite{Ebr} and \cite{Nkl3}.

In the light of Theorem \ref{ror} and Lemma \ref{fr}, it is reasonable to expect
that Einstein solvmanifolds are actually completely determined by their nilpotent
parts.  Let us describe this more precisely.

\begin{definition}
Given a metric nilpotent Lie algebra $(\ngo,\ip)$, a metric solvable Lie algebra
$(\sg=\ag\oplus\ngo,\ip')$ is called a {\it metric solvable extension} of
$(\ngo,\ip)$ if $\ngo$ is an ideal of $\sg$, $[\ag,\ag]\subset\ngo$ and
$\ip'|_{\ngo\times\ngo}=\ip$.
\end{definition}

It turns out that for each $(\ngo,\ip)$ there exists a unique rank-one metric
solvable extension of $(\ngo,\ip)$ which stands a chance of being an Einstein
manifold.  Indeed, since the inner product on $\ag$ is determined by Theorem
\ref{ror}, (iii), the only datum we need to recover $(\sg,\ip')$ from $(\ngo,\ip)$
is the way $\ag$ is acting on $\ngo$ by derivations.  Recall that $\dim{\ag}=1$.  It
follows from Lemma \ref{fr} that if $A\in\ag$ satisfies $||A||^2=\tr{\ad{A}}$ and we
set $D:=\ad{A}|_{\ngo}$, then the Ricci operator $\Ric_{\ip}$ of $(\ngo,\ip)$ equals
$$
\Ric_{\ip}=cI+D.
$$
But $\Ric_{\ip}$ is orthogonal to $D$ and actually to any symmetric derivation of
$(\ngo,\ip)$ (see \cite[(2)]{critical}), thus
\begin{equation}\label{pE}
\tr{DA}=-c\tr{A}, \qquad c=\tfrac{\tr{\Ric_{\ip}^2}}{\tr{\Ric_{\ip}}},
\end{equation}

\noindent for any symmetric derivation $A$ of $(\ngo,\ip)$.  This determines $D$ in
terms of $(\ngo,\ip)$ (recall that if $\ngo$ is nonabelian then $\tr{\Ric_{\ip}}<0$,
and in the abelian case both $\Ric_{\ip}$ and $D$ equal zero).  We have therefore
seen that a rank-one Einstein solvmanifold is completely determined by its (metric)
nilpotent part.  This fact turns the study of rank-one Einstein solvmanifolds into a
problem on nilpotent Lie algebras.

\begin{definition}\label{enil}
A nilpotent Lie algebra $\ngo$ is said to be an {\it Einstein nilradical} if it
admits an inner product $\ip$ such that there is a metric solvable extension of
$(\ngo,\ip)$ which is Einstein.
\end{definition}

Such a solvable extension must satisfy $[\ag,\ag]=0$ by Theorem \ref{stand} and it
can be assumed that $\dim{\ag}=1$ by Theorem \ref{ror}, (ii).  In other words,
Einstein nilradicals are precisely the {\it nilradicals} (i.e. the maximal nilpotent
ideal) of Lie algebras of Einstein solvmanifolds.

Given a nilpotent Lie algebra $\ngo$, an inner product $\ip$ on $\ngo$ can be
extended to construct an Einstein solvmanifold if and only if any of the following
equivalent conditions hold, which shows us that these left invariant metrics on
nilpotent Lie groups are very special from many other points of view:

\begin{itemize}
\item[ ]
\item[(i)] $\Ric_{\ip}=cI+D$ for some $c\in\RR$ and
$D\in\Der(\ngo)$, where $\Ric_{\ip}$ is the Ricci operator of $(\ngo,\ip)$ and
$\Der(\ngo)$ is the space of all derivations of $\ngo$.

\item[ ]
\item[(ii)] $\ip$ is a {\it Ricci soliton} metric: the solution $\ip_t$
with initial point $\ip_0=\ip$ to the Ricci flow
$$
\ddt\ip_t=-2\ricci_{\ip_t},
$$
remains isometric up to scaling to $\ip$, that is, $\ip_t=c_t\vp_t^*\ip$ for some
one parameter group of diffeomorphisms $\{\vp_t\}$ of $N$ and $c_t\in\RR$ (see for
instance \cite{libro}), where $N$ denotes the simply connected nilpotent Lie group
with Lie algebra $\ngo$.

\item[ ]
\item[(iii)] $\ip$ is a {\it quasi-Einstein} metric:
$$
\ricci_{\ip}=c\ip+L_X\ip
$$
for some $C^{\infty}$ vector field $X$ in $N$ and $c\in\RR$, where $L_X\ip$ denotes
the usual Lie derivative.  This class of metrics were actually first considered in
theoretical physics (see \cite{Frd} and \cite{ChvVln}).

\item[ ]
\item[(iv)] $\ip$ is a {\it minimal} metric:
$$
||\ricci_{\ip}||=\min \left\{ ||\ricci_{\ip'}|| :
\scalar(\ip')=\scalar(\ip)\right\},
$$
where $\ip'$ runs over all left invariant metrics on $N$ and $\scalar(\ip)$ denotes
scalar curvature of $(\ngo,\ip)$ (see \cite{praga}).  A nilpotent Lie group $N$ can
never admit an Einstein left invariant metric, unless it is abelian, and a way of
getting as close as possible to satisfy the Einstein condition is to have a minimal
metric. Indeed,
$$
||\ricci_{\ip}-\tfrac{\scalar(\ip)}{n}\ip||^2=||\ricci_{\ip}||^2
-\tfrac{\scalar(\ip)^2}{n}.
$$

\item[ ]
\end{itemize}

\begin{definition}\label{nilsoliton}
A left invariant metric on a nilpotent Lie group (or equivalently an inner product
on a nilpotent Lie algebra) is called a {\it nilsoliton} if it satisfies any of the
above conditions (i)-(iv).
\end{definition}

\begin{theorem}
A nilpotent Lie algebra $\ngo$ is an Einstein nilradical if and only if $\ngo$
admits a nilsoliton metric.
\end{theorem}

It is then reasonable to expect that uniqueness for nilsolitons should hold, as it
does for standard Einstein solvmanifolds (see Theorem \ref{u}).  This is actually
true, and even a very similar proof worked out.

\begin{theorem}\label{un}\cite{soliton}
There is at most one nilsoliton metric on a nilpotent Lie group up to isometry and
scaling.
\end{theorem}

We therefore obtain the following picture of one-to-one correspondences on the
classification problem:

$$
\begin{array}{c}
\hline \\ \\

\left\{\mbox{Rank-one Einstein solvmanifolds}\right\} / \mbox{isometry and scaling} \\ \\
\updownarrow \\ \\

\left\{\mbox{Ricci soliton (simply connected) nilmanifolds}\right\} / \mbox{isometry and scaling} \\ \\
\updownarrow \\ \\

\left\{\mbox{Einstein nilradicals}\right\} / \mbox{isomorphism} \\ \\
\hline \\
\end{array}
$$

Thus the classification of Einstein solvmanifolds reduces to a completely
`algebraic' problem; namely, the classification of nilpotent Lie algebras which are
Einstein nilradicals. This problem will be treated in Section \ref{clas}.

\section{Technical background}\label{Tb}

In this section, we fix the notation and give all the definitions and elementary
results we need to use throughout the paper.  All this mainly concerns the vector
space where the Lie algebras of a given dimension live and the action determining
the isomorphism relation between them.

Let us consider the space of all skew-symmetric algebras of dimension $n$, which is
parameterized by the vector space
$$
V=\lam=\{\mu:\RR^n\times\RR^n\longrightarrow\RR^n : \mu\; \mbox{bilinear and
skew-symmetric}\}.
$$
Then
$$
\nca=\{\mu\in V:\mu\;\mbox{satisfies Jacobi and is nilpotent}\}
$$
is an algebraic subset of $V$ as the Jacobi identity and the nilpotency condition
can both be written as zeroes of polynomial functions.  $\nca$ is often called the
{\it variety of nilpotent Lie algebras} (of dimension $n$).  There is a natural
action of $\G$ on $V$ given by
\begin{equation}\label{action}
g.\mu(X,Y)=g\mu(g^{-1}X,g^{-1}Y), \qquad X,Y\in\RR^n, \quad g\in\G,\quad \mu\in V.
\end{equation}

Recall that $\nca$ is $\G$-invariant and the Lie algebra isomorphism classes are
precisely the $\G$-orbits. The action of $\g$ on $V$ obtained by differentiation of
(\ref{action}) is given by
\begin{equation}\label{actiong}
\pi(\alpha)\mu=\alpha\mu(\cdot,\cdot)-\mu(\alpha\cdot,\cdot)-\mu(\cdot,\alpha\cdot),
\qquad \alpha\in\g,\quad\mu\in V.
\end{equation}

\noindent We note that $\pi(\alpha)\mu=0$ if and only if $\alpha\in\Der(\mu)$, the
Lie algebra of derivations of the algebra $\mu$, which is actually the Lie algebra
of $\Aut(\mu)$, the group of automorphisms of the algebra $\mu$.  Recall that
$\Aut(\mu)$ is the isotropy subgroup at $\mu$ for the action (\ref{action}), and so
$\dim{\G.\mu}=n^2-\dim{Der(\mu)}$.

The canonical inner product $\ip$ on $\RR^n$ determines an $\Or(n)$-invariant inner
product on $V$, also denoted by $\ip$, as follows:
\begin{equation}\label{innV}
\la\mu,\lambda\ra= \sum\limits_{ij}\la\mu(e_i,e_j),\lambda(e_i,e_j)\ra
=\sum\limits_{ijk}\la\mu(e_i,e_j),e_k\ra\la\lambda(e_i,e_j),e_k\ra,
\end{equation}

\noindent and also the standard $\Ad(\Or(n))$-invariant inner product on $\g$ given
by
\begin{equation}\label{inng}
\la \alpha,\beta\ra=\tr{\alpha \beta^{\mathrm t}}=\sum_{i}\la\alpha e_i,\beta e_i\ra
=\sum_{ij}\la\alpha e_i,e_j\ra\la\beta e_i,e_j\ra,, \qquad \alpha,\beta\in\g.
\end{equation}

\begin{remark} We have made several abuses of notation concerning inner products.  Recall
that $\ip$ has been used to denote an inner product on $\sg$, $\ngo$, $\RR^n$, $V$
and $\g$.
\end{remark}

We note that $\pi(\alpha)^t=\pi(\alpha^t)$ and $(\ad{\alpha})^t=\ad{\alpha^t}$ for
any $\alpha\in\g$, due to the choice of canonical inner products everywhere.

Let $\tg$ denote the set of all diagonal $n\times n$ matrices.  If $\{
e_1',...,e_n'\}$ is the basis of $(\RR^n)^*$ dual to the canonical basis $\{
e_1,...,e_n\}$ then
$$
\{ v_{ijk}=(e_i'\wedge e_j')\otimes e_k : 1\leq i<j\leq n, \; 1\leq k\leq n\}
$$
is a basis of weight vectors of $V$ for the action (\ref{action}), where $v_{ijk}$
is actually the bilinear form on $\RR^n$ defined by
$v_{ijk}(e_i,e_j)=-v_{ijk}(e_j,e_i)=e_k$ and zero otherwise.  The corresponding
weights $\alpha_{ij}^k\in\tg$, $i<j$, are given by
\begin{equation}\label{alfas}
\pi(\alpha)v_{ijk}=(a_k-a_i-a_j)v_{ijk}=\la\alpha,\alpha_{ij}^k\ra v_{ijk},
\quad\forall\alpha=\left[\begin{smallmatrix} a_1&&\\ &\ddots&\\ &&a_n
\end{smallmatrix}\right]\in\tg,
\end{equation}

\noindent where $\alpha_{ij}^k=E_{kk}-E_{ii}-E_{jj}$ and $\ip$ is the inner product
defined in (\ref{inng}).  As usual $E_{rs}$ denotes the matrix whose only nonzero
coefficient is $1$ at entry $rs$.  The structural constants $\mu_{ij}^k$ of an
algebra $\mu\in V$ are then given by $$ \mu(e_i,e_j)=\sum_{k=1}^n\mu_{ij}^ke_k,
\qquad\mbox{or}\quad \mu=\sum\mu_{ij}^kv_{ijk}, \quad i<j.
$$
Let $\tg^+$ denote the Weyl chamber of $\g$ given by
\begin{equation}\label{weyl}
\tg^+=\left\{\left[\begin{smallmatrix} a_1&&\\ &\ddots&\\ &&a_n
\end{smallmatrix}\right]\in\tg:a_1\leq...\leq a_n\right\}.
\end{equation}

For $\alpha\in\tg^+$ we define the parabolic subgroup
\begin{equation}\label{para}
P_{\alpha}=B\G_{\alpha},
\end{equation}

\noindent where $B$ is the subgroup of $\G$ of all lower triangular invertible
matrices and
$$
\G_{\alpha}=\{ g\in\G:g\alpha g^{-1}=\alpha\}.
$$
In general, for any $\alpha'\in\g$ which is diagonalizable over $\RR$, we let
$P_{\alpha'}:=gP_{\alpha}g^{-1}$ if $\alpha'=g\alpha g^{-1}$, $\alpha\in\tg^+$. This
is well defined since $h\alpha h^{-1}=g\alpha g^{-1}$ implies that
$h^{-1}g\in\G_{\alpha}\subset P_{\alpha}$ and so
$h^{-1}gP_{\alpha}g^{-1}h=P_{\alpha}$.

There is an ordered basis of $V$ with respect to which the action of $g$ on $V$ is
lower triangular for any $g\in B$, and furthermore the eigenvalues of $\pi(\alpha)$
are increasing for any $\alpha\in\tg^+$.

Given a finite subset $X$ of $\tg$, we denote by $\CH(X)$ the convex hull of $X$ and
by $\mcc(X)$ the {\it minimal convex combination of} $X$, that is, the (unique)
vector of minimal norm (or closest to the origin) in $\CH(X)$.  If
$X=\{\alpha_1,...,\alpha_r\}\subset\tg$ and $\beta:=\mcc(X)$, then there exist
$c_i\geq 0$, $i=1,...,r$, such that $\sum\limits_{i=1}^rc_i=1$ and
$\beta=\sum\limits_{i=1}^rc_i\alpha_i$. Since $\la\beta,\alpha_i\ra\geq ||\beta||^2$
for all $i$ (why?), we have that
$$
||\beta||^2=\sum_{i=1}^rc_i\la\beta,\alpha_i\ra\geq ||\beta||^2,
$$
from which follows that $\la\beta,\alpha_i\ra=||\beta||^2$ for all $i$ such that $c_i>0$.  We
can therefore assume that $\la\beta,\alpha_i\ra=||\beta||^2$ for all $i$, and also
that $\beta=\mcc(\{\alpha_1,...,\alpha_s\})$, where $\{\alpha_1,...,\alpha_s\}$ is a
linearly independent subset of $X$.  Thus the $s\times s$ matrix
$U:=\left[\la\alpha_i,\alpha_j\ra\right]$ is invertible and satisfies
\begin{equation}\label{rat}
U\left[\begin{smallmatrix} c_1 \\ \vdots \\ c_s\end{smallmatrix}\right]=
||\beta||^2\left[\begin{smallmatrix} 1 \\ \vdots \\ 1\end{smallmatrix}\right].
\end{equation}

In particular, if all the entries of $\alpha_i$ are in $\QQ$ for any $i=1,...,r$,
then also the entries of $\beta:=\mcc(X)$ are all in $\QQ$.  Indeed,
$\tfrac{c_i}{||\beta||^2}\in\QQ$ for all $i$ and so their sum
$\tfrac{1}{||\beta||^2}\in\QQ$, which implies that $c_i\in\QQ$ for all $i$ and
consequently $\beta$ has all its coefficients in $\QQ$.

\section{Variational approach to Einstein solvmanifolds}\label{va}

Einstein metrics are often considered as the nicest, or most privileged ones on a
given differentiable manifold (see for instance \cite[Introduction]{Bss}).  One of
the justifications is the following result due to Hilbert (see \cite{Hlb}): the
Einstein condition for a compact Riemannian manifold $(M,g_{\circ})$ of volume one
is equivalent to the fact that the total scalar curvature functional
$$
\scalar : g\mapsto\int_M\scalar(g)\mu_g
$$
admits $g_{\circ}$ as a critical point on the space of all metrics of volume one
(see also \cite[4.21]{Bss}).  This variational approach still works for $G$-invariant
metrics on $M$, where $G$ is any compact Lie group acting transitively on $M$ (see
\cite[4.23]{Bss}).

On the other hand, it is proved in \cite{Jns} that in a unimodular $n$-dimensional
Lie group, the Einstein left invariant metrics are precisely the critical points of
the scalar curvature functional on the set of all left invariant metrics having a
fixed volume element.  However, this fails in the non-unimodular case.  For
instance, if $\sg$ is a solvable non-unimodular Lie algebra, then the scalar
curvature functional restricted to any leaf
$F=\{t\}\times\Sl(\sg)/\SO(\sg)\subset\pca$ of inner products, has no critical
points (see \cite[3.5]{Hbr}).  Thus, the approach to study Einstein solvmanifolds by
a variational method should be different.

In this section, we shall describe the approach proposed in the introduction: to
vary Lie brackets rather than inner products.  Recall that when $\ngo$ is an
$n$-dimensional nilpotent Lie algebra, then the set of all inner products on $\ngo$
is very nice, it is parameterized by the symmetric space $\Gl_n(\RR)/\Or(n)$.
However, isometry classes are precisely the orbits of the action on
$\Gl_n(\RR)/\Or(n)$ of the group of automorphisms $\Aut(\ngo)$, a group mostly
unknown, hard to compute, and far from being reductive, that is, ugly from the point
of view of invariant theory.  If we instead vary Lie brackets, isometry classes will
be given by $\Or(n)$-orbits, a beautiful group.  But since nothing is for free in
mathematics, the set of left invariant metrics will now be parameterized by a
$\Gl_n(\RR)$-orbit in the variety $\nca$ of $n$-dimensional nilpotent Lie algebras,
a terrible space.

We fix an inner product vector space
$$
(\sg=\RR H\oplus\RR^n,\ip),\qquad \la H,\RR^n\ra=0,\quad \la H,H\ra=1,
$$
such that the restriction $\ip|_{\RR^n\times\RR^n}$ is the canonical inner product
on $\RR^n$, which will also be denoted by $\ip$.  A linear operator on $\RR^n$ will
be sometimes identified with its matrix in the canonical basis $\{ e_1,...,e_n\}$ of
$\RR^n$.  The metric Lie algebra corresponding to any $(n+1)$-dimensional rank-one
solvmanifold, can be modeled on $(\sg=\RR H\oplus\ngo,\ip)$ for some nilpotent Lie
bracket $\mu$ on $\RR^n$ and some $D\in\Der(\mu)$, the space of derivations of
$(\RR^n,\mu)$.  Indeed, these data define a solvable Lie bracket $[\cdot,\cdot]$ on
$\sg$ by
\begin{equation}\label{solv}
[H,X]=DX,\qquad [X,Y]=\mu(X,Y), \qquad X,Y\in\RR^n,
\end{equation}

\noindent and the solvmanifold is then the simply connected Lie group $S$ with Lie
algebra $(\sg,[\cdot,\cdot])$ endowed with the left invariant Riemannian metric
determined by $\ip$.  We shall assume from now on that $\mu\ne 0$ since the case
$\mu=0$ (i.e. abelian nilradical) is well understood (see \cite[Proposition
6.12]{Hbr}).  We have seen in the paragraph above Definition \ref{enil} that for a
given $\mu$, there exists a unique  symmetric derivation $D_{\mu}$ to consider if we
want to get Einstein solvmanifolds.  We can therefore associate to each nilpotent
Lie bracket $\mu$ on $\RR^n$ a distinguished rank-one solvmanifold $S_{\mu}$,
defined by the data $\mu,D_{\mu}$ as in (\ref{solv}), which is the only one with a
chance of being Einstein among all those metric solvable extensions of $(\mu,\ip)$.

We note that conversely, any $(n+1)$-dimensional rank-one Einstein solvmanifold is
isometric to $S_{\mu}$ for some nilpotent $\mu$. Thus the set $\nca$ of all
nilpotent Lie brackets on $\RR^n$ parameterizes a space of $(n+1)$-dimensional
rank-one solvmanifolds
$$
\{ S_{\mu}:\mu\in\nca\},
$$
containing all those which are Einstein in that dimension.

Concerning the identification
$$
\mu\longleftrightarrow (N_{\mu},\ip),
$$
where $N_{\mu}$ is the simply connected nilpotent Lie group with Lie algebra
$(\RR^n,\mu)$, the $\G$-action on $\nca$ defined in (\ref{action}) has the following
geometric interpretation: each $g\in\G$ determines a Riemannian isometry
\begin{equation}\label{id}
(N_{g.\mu},\ip)\longrightarrow (N_{\mu},\la g\cdot,g\cdot\ra)
\end{equation}

\noindent by exponentiating the Lie algebra isomorphism
$g^{-1}:(\RR^n,g.\mu)\longrightarrow(\RR^n,\mu)$.  Thus the orbit $\G.\mu$ may be
viewed as a parametrization of the set of all left invariant metrics on $N_{\mu}$.
By a result of E. Wilson, two pairs $(N_{\mu},\ip)$, $(N_{\lambda},\ip)$ are
isometric if and only if $\mu$ and $\lambda$ are in the same $\Or(n)$-orbit (see
\cite[Appendix]{minimal}), where $\Or(n)$ denotes the subgroup of $\G$ of orthogonal
matrices. Also, two solvmanifolds $S_{\mu}$ and $S_{\lambda}$ with
$\mu,\lambda\in\nca$ are isometric if and only if there exists $g\in\Or(n)$ such
that $g.\mu=\lambda$ (see \cite[Proposition 4]{critical}).  From (\ref{id}) and the
definition of $S_{\mu}$ we obtain the following result.

\begin{lemma}\label{enilmu}
If $\mu\in\nca$ then the nilpotent Lie algebra $(\RR^n,\mu)$ is an Einstein
nilradical if and only if $S_{g.\mu}$ is Einstein for some $g\in\G$.
\end{lemma}

Recall that being an Einstein nilradical is a property of a whole $\G$-orbit in
$\nca$, that is, of the isomorphism class of a given $\mu$.

For any $\mu\in\nca$ we have that the scalar curvature of $(N_{\mu},\ip)$ is given
by $\scalar(\mu)=-\unc ||\mu||^2$, which says that normalizing by scalar curvature
and by the spheres of $V$ is actually equivalent.  The critical points of any
scaling invariant curvature functional on $\nca$ appear then as very natural
candidates to be distinguished left invariant metrics on nilpotent Lie groups.

\begin{theorem}\label{crit}\cite{soliton,critical,einsteinsolv}
For a nonzero $\mu\in\nca$, the following conditions are equivalent:
\begin{itemize}
\item[(i)] $S_{\mu}$ is Einstein.

\item[(ii)] $(N_{\mu},\ip)$ is a nilsoliton.

\item[(iii)] $\mu$ is a critical point of the functional $F:V\longrightarrow\RR$ defined by
$$
F(\mu)=\tfrac{16}{||\mu||^4}\tr{\Ric_{\mu}^2},
$$
where $\Ric_{\mu}$ denotes the Ricci operator of $(N_{\mu},\ip)$.

\item[(iv)] $\mu$ is a minimum of $F|_{\G.\mu}$ (i.e. $(N_{\mu},\ip)$ is minimal).

\item[(v)] $\Ric_{\mu}\in\RR I\oplus\Der(\mu)$.
\end{itemize}
Under these conditions, the set of critical points of $F$ lying in $\G.\mu$ equals
$\Or(n).\mu$ (up to scaling).
\end{theorem}

Thus another natural approach to find rank-one Einstein solvmanifolds would be to
use the negative gradient flow of the functional $F$.  It follows from \cite[Lemma
6]{critical} that if $\pi$ is the representation defined in (\ref{actiong}) then
$$
\grad(F)_{\mu}=\tfrac{16}{||\mu||^6}\left(||\mu||^2\pi(\Ric_{\mu})\mu-4\tr{\Ric_{\mu}^2}\mu\right).
$$
Since $F$ is invariant under scaling we know that $||\mu||$ will remain constant in
time along the flow. We may therefore restrict ourselves to the sphere of radius
$2$, where the negative gradient flow $\mu=\mu(t)$ of $F$ becomes
\begin{equation}\label{flow}
    \ddt\mu=-\pi(\Ric_{\mu})\mu+\tr{\Ric_{\mu}^2}\mu.
\end{equation}

Notice that $\mu(t)$ is a solution to this differential equation if and only if
$g.\mu(t)$ is so for any $g\in\Or(n)$, according to the $\Or(n)$-invariance of $F$.
The existence of $\lim\limits_{t\to\infty}\mu(t)$ is guaranteed by the compactness
of the sphere and the fact that $F$ is a polynomial (see for instance \cite[Section
2.5]{Sjm}).

\begin{lemma}\label{strataflow}\cite{einsteinsolv}
For $\mu_0\in V$, $||\mu_0||=2$, let $\mu(t)$ be the flow defined in {\rm
(\ref{flow})} with $\mu(0)=\mu_0$ and put $\lambda=\lim\limits_{t\to\infty}\mu(t)$.
Then
\begin{itemize}
    \item[(i)] $\mu(t)\in\G.\mu_0$ for all $t$.
    \item[(ii)] $\lambda\in\overline{\G.\mu_0}$.
    \item[(iii)] $S_{\lambda}$ is Einstein.
   \end{itemize}
\end{lemma}

Part (i) follows from the fact that $\ddt\mu\in\tang_{\mu}\G.\mu$ for all $t$ (see
(\ref{flow})), and part (ii) is just a consequence of (i).  Condition (ii) is often
referred in the literature as the Lie algebra $\mu_0$ {\it degenerates} to the Lie
algebra $\lambda$.  Some interplays between degenerations and Riemannian geometry of
Lie groups have been explored in \cite{inter}, by using the fact that for us, the
orbit $\G.\mu_0$ is the set of all left invariant metrics on $N_{\mu_0}$.  We note
that if the limit $\lambda\in\G.\mu_0$, then $\mu_0$ is an Einstein nilradical.  We
do not know if the converse holds.  Since $\lambda$ is a critical point of $F$ and
$\lambda\in\nca$ by (ii) and the fact that $\nca$ is closed, we have that part (iii)
follows from Theorem \ref{crit}.

In geometric invariant theory, a moment map for linear reductive Lie group actions
over $\CC$ has been defined in \cite{Nss} and \cite{Krw1} (see Appendix).  In our
situation, it is an $\Or(n)$-equivariant map
$$
m:V\smallsetminus\{ 0\}\longrightarrow\sym(n),
$$
defined implicitly by
\begin{equation}\label{defmm}
\la m(\mu),\alpha\ra=\tfrac{1}{||\mu||^2}\la\pi(\alpha)\mu,\mu\ra, \qquad \mu\in
V\smallsetminus\{ 0\}, \; \alpha\in\sym(n).
\end{equation}

\noindent We are using $\g=\sog(n)\oplus\sym(n)$ as the Cartan decomposition for the
Lie algebra $\g$ of $\G$, where $\sog(n)$ and $\sym(n)$ denote the subspaces of
skew-symmetric and symmetric matrices, respectively.

Recall that $\nca\subset V$ and each $\mu\in\nca$ determines two Riemannian
manifolds $S_{\mu}$ and $(N_{\mu},\ip)$. A remarkable fact is that this moment map
encodes geometric information on $S_{\mu}$ and $(N_{\mu},\ip)$; indeed, it was
proved in \cite{minimal} that
\begin{equation}\label{mmR}
m(\mu)=\tfrac{4}{||\mu||^2}\Ric_{\mu}.
\end{equation}

This allows us to use strong and well-known results on the moment given in
\cite{Krw1} and \cite{Nss}, and proved in \cite{Mrn} for the real case (see the
Appendix for an overview on such results).  We note that the functional $F$ defined
in Theorem \ref{crit}, (iii) is precisely $F(\mu)=||m(\mu)||^2$, and so the
equivalence between (iii) and (iv) in Theorem \ref{crit} follows from Theorem
\ref{marian}, (i).  It should be pointed out that actually most of the results in
Theorem \ref{crit} follow from general results on the moment map proved in
\cite{Mrn}.  For instance, the last sentence about uniqueness of critical points of
$F$ (see Theorem \ref{marian}, (ii)), is easily seen to be equivalent to the
uniqueness of standard Einstein solvmanifolds (see Theorem \ref{u}) and nilsolitons
(see Theorem \ref{un}).

In Section \ref{st}, we shall see that one can go further in the application of
geometric invariant theory to the study of Einstein solvmanifolds, by considering a
stratification for $\nca$ intimately related to the moment map.

\section{On the classification of Einstein solvmanifolds}\label{clas}

As we have seen in Section \ref{pre}, the classification of Einstein solvmanifolds
is essentially reduced to the rank-one case.  There is a bijection between the set
of all isometry classes of rank-one Einstein solvmanifolds and the set of isometry
classes of certain distinguished left invariant metrics on nilpotent Lie groups
called nilsolitons, and the uniqueness up to isometry of nilsolitons finally
determines a new bijection with the set of all isomorphism classes of Einstein
nilradicals.  For better or worse, what we get in the end is then a classification
problem on nilpotent Lie algebras.

Recall that a nilpotent Lie algebra $\ngo$ is an Einstein nilradical if and only if
$\ngo$ admits a nilsoliton, that is, an inner product $\ip$ such that the
corresponding Ricci operator $\Ric_{\ip}$ satisfies
$$
\Ric_{\ip}=cI+D, \qquad \mbox{for some}\; c\in\RR,\; D\in\Der(\ngo).
$$
Therefore, in order to understand or classify Einstein nilradicals, a main problem
would be how to translate this condition based on the existence of an inner product
on $\ngo$ having a certain property into purely Lie theoretic conditions on $\ngo$.
The following questions also arise:

\begin{itemize}
\item[(A)] Besides the existence of an $\NN$-gradation, is there any other neat structural obstruction for a nilpotent Lie algebra to be an Einstein nilradical?

\item[ ]\item[(B)] Is there any algebraic condition on a nilpotent Lie algebra which is sufficient to be an Einstein nilradical?

\item[ ]\item[(C)] An $\NN$-graded nilpotent Lie algebra can or can not be an
Einstein nilradical, what is more likely?
\end{itemize}

Let us now review what we do know on the classification of Einstein nilradicals.

Any nilpotent Lie algebra of dimension $\leq 6$ is an Einstein nilradical (see
\cite{Wll}). There are $34$ of them in dimension $6$,  giving rise to $29$ different
eigenvalue-types (there are $5$ eigenvalue-types with exactly two algebras).  In
dimension $7$, the first nilpotent Lie algebras without any $\NN$-gradation appear,
but also do the first examples of $\NN$-graded Lie algebras which are not Einstein
nilradicals.  The family of $7$-dimensional nilpotent Lie algebras defined for any
$t\in\RR$ by
\begin{equation}\label{ex7}
\begin{array}{lll}
[X_1,X_2]_t=X_3, & [X_1,X_5]_t=X_6, & [X_2,X_4]_t=X_6, \\ \\

[X_1,X_3]_t=X_4, & [X_1,X_6]_t=X_7, &  [X_2,X_5]_t=tX_7, \\ \\

[X_1,X_4]_t=X_5, & [X_2,X_3]_t=X_5, &  [X_3,X_4]_t=(1-t)X_5,
\end{array}
\end{equation}

\noindent is really a curve in the set of isomorphism classes of algebras (i.e.
$\lb_t\simeq\lb_s$ if and only if $t=s$) and $\lb_t$ turns to be an Einstein
nilradical if and only if $t\ne 0,1$ (see \cite{einsteinsolv}).  Recall that all of
them admit the gradation $\ngo=\ngo_1\oplus\ngo_2\oplus...\oplus\ngo_7$, $\ngo_i=\RR
X_i$ for all $i$.  This example in particular shows that to be an Einstein
nilradical is not a property which depends continuously on the structural constants
of the Lie algebra.

Perhaps the nicest source of examples of Einstein nilradicals is the following.

\begin{theorem}\cite{Tmr}
Let $\ggo$ be a real semisimple Lie algebra.  Then the nilradical of any parabolic
subalgebra of $\ggo$ is an Einstein nilradical.
\end{theorem}

If we add to this that H-type Lie algebras and any nilpotent Lie algebra admitting a
naturally reductive left invariant metric are Einstein nilradicals, one may get the
impression that any nilpotent Lie algebra which is special or distinguished in some
way, or just has a `name', will be an Einstein nilradical.  This is contradicted by
the following surprising result, which asserts that free nilpotent Lie algebras are
rarely Einstein nilradicals.

\begin{theorem}\cite{Nkl1}
A free $p$-step nilpotent Lie algebra on $m$ generators is an Einstein nilradical if
and only if
\begin{itemize}
\item $p=1,2$;
\item $p=3$ and $m=2,3,4,5$;
\item $p=4$ and $m=2$;
\item $p=5$ and $m=2$.
\end{itemize}
\end{theorem}

A nilpotent Lie algebra $\ngo$ is said to be {\it filiform} if $\dim{\ngo}=n$ and
$\ngo$ is $(n-1)$-step nilpotent.  These algebras may be seen as those which are as
far as possible from being abelian along the class of nilpotent Lie algebras, and in
fact most of them admit at most one $\NN$-gradation.  Several families of filiform
algebras which are not Einstein nilradicals have been found in \cite{Nkl2}, as well
as many isolated examples of non-Einstein nilradicals belonging to a curve of
Einstein nilradicals as in example (\ref{ex7}).  In \cite{Arr}, a weaker version of
Theorem \ref{Upos} given in \cite{Nkl2} is used to get a classification of
$8$-dimensional filiform Einstein nilradicals.

The lack of $\NN$-gradations is not however the only obstacle one can find for
Einstein nilradicals.  Several examples of non-Einstein nilradicals are already
known in the class of $2$-step nilpotent Lie algebras (i.e. $[\ngo,[\ngo,\ngo]]=0$),
the closest ones to being abelian and so algebras which usually admit plenty of
different $\NN$-gradations.

\begin{definition} A $2$-step nilpotent Lie algebra $\ngo$ is said to be of {\it type} $(p,q)$ if $\dim{\ngo}=p+q$ and $\dim{[\ngo,\ngo]}=p$.
\end{definition}

In \cite{einsteinsolv}, certain $2$-step nilpotent Lie algebras attached to graphs
are considered (of type $(p,q)$ if the graph has $q$ vertices and $p$ edges) and it
is proved that they are Einstein nilradicals if and only if the graph is positive
(i.e. when certain uniquely defined weighting on the set of edges is positive).  For
instance, any regular graph and also any tree such that any of its edges is adjacent
to at most three other edges is positive.  On the other hand, a graph is not
positive under the following condition: there are two joined vertices $v$ and $w$
such that $v$ is joined to $r$ vertices of valency $1$, $w$ is joined to $s$
vertices of valency $1$, both are joined to $t$ vertices of valency $2$ and
$(r,s,t)$ is not in a set of only a few exceptional small triples. This provides a
great deal of $2$-step non-Einstein nilradicals, starting from types $(5,6)$ and
$(7,5)$, and any dimension $\geq 11$  is attained.

Many other $2$-step algebras of type $(6,5)$ and $(7,5)$ which are not Einstein
nilradicals have appeared from the complete classification for types $(p,q)$ with
$q\leq 5$ and $(p,q)\ne (5,5)$ carried out in \cite{Nkl3}.

Curiously enough, at this point of the story, with so many examples of non-Einstein
nilradicals available, a curve was still missing.  In each fixed dimension,
only finitely many nilpotent Lie algebras which are not Einstein nilradicals have
showed up.  But this potential candidate to a conjecture has recently been dismissed
by the following result.

\begin{theorem}\cite{Wll2}
Let $\ngo_t$ be the $9$-dimensional Lie algebra with Lie bracket defined by
$$
\begin{array}{lll}

[X_5,X_4]_t=X_7, & [X_1,X_6]_t=X_8, &  [X_3,X_2]_t=X_9, \\ \\

[X_3,X_6]_t=tX_7, & [X_5,X_2]_t=tX_8, &  [X_1,X_4]_t=tX_9, \\ \\

[X_1,X_2]_t=X_7. &&
\end{array}
$$

Then $\ngo_{t},$  $t\in (1,\infty)$, is a curve of pairwise non-isomorphic $2$-step
nilpotent Lie algebras of type $(3,6)$, none of which is an Einstein nilradical.
\end{theorem}

The following definition is motivated by (\ref{pE}), a condition a rank-one solvable
extension of a nilpotent Lie algebra must satisfy in order to have a chance of being
Einstein.

\begin{definition}\label{preE} A derivation $\phi$ of a real Lie algebra $\ggo$ is called {\it pre-Einstein} if it is diagonalizable over $\RR$ and
$$
\tr{\phi\psi}=\tr{\psi}, \qquad\forall\; \psi\in\Der(\ggo).
$$
\end{definition}

The following result is based on the fact that $\Aut(\ggo)$ is an algebraic group.

\begin{theorem}\cite{Nkl3}\label{eupreE}
Any Lie algebra $\ggo$ admits a pre-Einstein derivation, which is unique up to
$\Aut(\ggo)$-conjugation and has eigenvalues in $\QQ$.
\end{theorem}

Let $\ngo$ be a nilpotent Lie algebra with pre-Einstein derivation $\phi$.  We note
that if $\ngo$ admits a nilsoliton metric, say with $\Ric_{\ip}=cI+D$, then $D$
necessarily equals $\phi$ up to scaling and conjugation (see (\ref{pE})), and thus
the eigenvalue-type of the corresponding Einstein solvmanifold is the set of
eigenvalues of $\phi$ up to scaling.  In particular, $\phi>0$.  It is proved in
\cite{Nkl3} that also $\ad{\phi}\geq 0$ as long as $\ngo$ is an Einstein nilradical.
These conditions are not however sufficient to guarantee that $\ngo$ is an Einstein
nilradical (see \cite{Nkl0}).  In order to get a necessary and sufficient condition
in terms of $\phi$ we have to work harder.

Let us first consider
\begin{equation}\label{gphi}
\ggo_{\phi}:=\{\alpha\in\glg(\ngo):[\alpha,\phi]=0, \quad\tr{\alpha\phi}=0,
\quad\tr{\alpha}=0\}
\end{equation}

\noindent and let $G_{\phi}$ be the connected Lie subgroup of $\Gl(\ngo)$ with Lie
algebra $\ggo_{\phi}$.  Recall that the Lie bracket $\lb$ of $\ngo$ belongs to the
vector space $\Lambda^2\ngo^*\otimes\ngo$ of skew-symmetric bilinear maps from
$\ngo\times\ngo$ to $\ngo$, on which $\Gl(\ngo)$ is acting naturally by
$g.\lb=g[g^{-1}\cdot,g^{-1}\cdot]$.

\begin{theorem}\cite{Nkl3}\label{madreN}
Let $\ngo$ be a nilpotent Lie algebra with pre-Einstein derivation $\phi$. Then
$\ngo$ is an Einstein nilradical if and only if the orbit $G_{\phi}.\lb$ is closed
in $\Lambda^2\ngo^*\otimes\ngo$.
\end{theorem}

This is certainly the strongest general result we know so far concerning questions
(A) and (B) above, and of course it has many useful applications, some of which we will now
describe (see also Theorem \ref{madre} for a turned to be equivalent result).

\begin{definition}\label{nice}  Let $\{ X_1,...,X_n\}$ be a basis for a nilpotent Lie algebra $\ngo$,
with structural constants $c_{ij}^k$'s given by
$[X_i,X_j]=\sum\limits_{k=1}^{n}c_{ij}^kX_k$.  Then the basis $\{ X_i\}$ is said to
be {\it nice} if the following conditions hold:
\begin{itemize}
    \item for all $i<j$ there is at most one $k$ such that $c_{ij}^k\ne 0$,
    \item if $c_{ij}^k$ and $c_{i'j'}^k$ are nonzero then either $\{ i,j\}=\{
    i',j'\}$ or $\{ i,j\}\cap\{ i',j'\}=\emptyset$.
\end{itemize}
\end{definition}

A nice property a nice basis $\{ X_i\}$ has is that the Ricci operator $\Ric_{\ip}$
of any inner product $\ip$ for which $\{ X_i\}$ is orthogonal diagonalizes with
respect to $\{ X_i\}$ (see \cite[Lemma 3.9]{einsteinsolv}).  Uniform bases
considered in \cite{Dlf} and \cite{Wlt} are nice.  The existence of a nice basis for
a nilpotent Lie algebra looks like a strong condition, although we do not know of
any example for which we can prove the non-existence of a nice basis.  Not even an
existence result for such example is available.

Let $(\ngo,\ip)$ be a metric nilpotent Lie algebra with orthogonal basis $\{ X_i\}$
and structural constants $[X_i,X_j]=\sum\limits_{k=1}^{n}c_{ij}^kX_k$.  If we fix an
enumeration of the set $\left\{\alpha_{ij}^k:c_{ij}^k\ne 0\right\}$ (see Section
\ref{Tb}), we can define the symmetric matrix
\begin{equation}\label{defU}
U=\left[\left\langle\alpha_{ij}^k,\alpha_{i'j'}^{k'}\right\rangle\right],
\end{equation}

\noindent and state the following useful result.

\begin{theorem}\cite{Pyn}\label{tracy}
Assume that $(\ngo,\ip)$ satisfies $\Ric_{\ip}\in\tg$.  Then $(\ngo,\ip)$ is a
nilsoliton if and only if
$$
U\left[(c_{ij}^k)^2\right]=c [1], \qquad c\in\RR,
$$
where $\left[(c_{ij}^k)^2\right]$ is meant as a column vector in the same order used
in {\rm (\ref{defU})} for defining $U$ and $[1]$ is the column vector with all
entries equal to $1$.
\end{theorem}

It turns out that equations $U\left[(c_{ij}^k)^2\right]=c[1]$ are precisely those given by the
Lagrange method applied to find critical points of the functional $F$ in Theorem
\ref{crit}.  In \cite{Pyn}, a Cartan matrix is associated to $U$ and the theory of
Kac-Moody algebras is applied to analyze the solutions space of such a linear
system.

Recall that $\mcc(X)$ denotes the unique vector of minimal norm in the convex hull
$\CH(X)$ of a finite subset $X$ of $\tg$.

\begin{theorem}\cite{Nkl3}\label{Upos}
A nonabelian nilpotent Lie algebra $\ngo$ with a nice basis $\{ X_i\}$ and
structural constants $[X_i,X_j]=\sum\limits_{k=1}^{n}c_{ij}^kX_k$ is an Einstein
nilradical if and only if any of the following equivalent conditions hold:
\begin{itemize}
\item[(i)] $\mcc\{\alpha_{ij}^k:c_{ij}^k\ne 0\}$ lies in the interior of $\CH\left(\{\alpha_{ij}^k:c_{ij}^k\ne 0\}\right)$.

\item[(ii)] Equation $U[x_{ij}^k]=[1]$ has a positive solution $[x_{ij}^k]$.
\end{itemize}
\end{theorem}

This is a non-constructive result, in the sense that it is in general very difficult
to explicitly find the nilsoliton metric.  The absence of an inner product in its
statement (compare with Theorem \ref{tracy}), however, makes of Theorem \ref{Upos},
quite a useful result.

\begin{theorem}\cite{Nkl3}
Let $\ngo_1$, $\ngo_2$ be real nilpotent Lie algebras which are isomorphic as
complex Lie algebras (i.e. they have isomorphic complexifications
$\ngo_{i}\otimes\CC$).  Then $\ngo_1$ is an Einstein nilradical if and only if
$\ngo_2$ is so, and in that case, they have the same eigenvalue-type.
\end{theorem}

This turns our classification of Einstein nilradicals into a problem on complex
nilpotent Lie algebras, with all the advantages an algebraically closed field has if
we want to use known classifications in the literature or results from algebraic
geometry and geometric invariant theory.

The following result reduces the classification of Einstein nilradicals to those
which are {\it indecomposable} (i.e. non-isomorphic to a direct sum of Lie
algebras).

\begin{theorem}\cite{Nkl3}\label{suma}
Let $\ngo=\ngo_1\oplus\ngo_2$ be a nilpotent Lie algebra which is the direct sum of
two ideals $\ngo_1$ and $\ngo_2$. Then $\ngo$ is an Einstein nilradical if and only
if both $\ngo_1$ and $\ngo_2$ are Einstein nilradicals.
\end{theorem}

Any $2$-step nilpotent Lie algebra of type $(p,q)$ can be identified with an element
in the vector space $V_{q,p}:=\Lambda^2(\RR^q)^*\otimes\RR^p$, and it is easy to see
that two of them are isomorphic if and only if they lie in the same
$\Gl_q(\RR)\times\Gl_p(\RR)$-orbit.

\begin{theorem}\label{generic}\cite{Ebr, Nkl3}
If $(p,q)\ne (2,2k+1)$, then the
vector space $V_{q,p}$ contains an open and dense subset of Einstein nilradicals of
eigenvalue-type $(1<2;q,p)$.
\end{theorem}

In view of this result, one may say that for most types, a $2$-step nilpotent Lie
algebra of type $(p,q)$ is typically, or generically, an Einstein nilradical of
eigenvalue-type $(1<2;q,p)$.  This is no doubt an important indicator related to
question (C) above, but we must go carefully. What Theorem \ref{generic} is actually
asserting is that if one throws a dart on $V_{q,p}$, then, with probability one, the
dart will hit at a Lie bracket $\lb\in V_{q,p}$ which is an Einstein nilradical of
eigenvalue-type $(1<2;q,p)$.  Recall that each algebra of type $(p,q)$ is identified
with a whole $\Gl_q(\RR)\times\Gl_p(\RR)$-orbit in $V_{q,p}$, not with a single
point, and some of these orbits can be much thicker than others.

Let us consider a simple example to illustrate this phenomenon.  There are exactly
$7$ algebras up to isomorphism in the vector space $V_{4,2}$, including the abelian
one.  Only two of them are Einstein nilradicals of type $(1<2;4,2)$; namely, the
H-type Lie algebra $\hg_3\otimes\CC$ (i.e. the complexification of $\hg_3$ viewed as
real) and $\hg_3\oplus\hg_3$, where $\hg_3$ denotes the $3$-dimensional Heisenberg
algebra (see \cite[Table 4]{Wll}).  If we fix basis $\{ X_1,...,X_4\}$ and $\{
Z_1,Z_2\}$ of $\RR^4$ and $\RR^2$, respectively, then each $\lb\in V_{4,2}$ is
determined by $12$ structural constants as follows:
$$
[X_i,X_j]=c_{ij}^1Z_1+c_{ij}^2Z_2, \qquad c_{ij}^k\in\RR, \quad 1\leq i<j\leq 4,
\quad k=1,2.
$$
If we take variables $x,y$ and define the skew-symmetric matrix $J$ with $ij$ entry,
$i<j$, given by $c_{ij}^1x+c_{ij}^2y$, then $\det{J}$ is a $4$-degree homogeneous
polynomial on $(x,y)$ with a `square root' $f(x,y)$, a $2$-degree homogeneous
polynomial called the Pfaffian form of $\lb$ (see \cite[Section 2]{ratform}).  Thus
the Hessian of $f$ is a real number $h(\lb)$ which depends polynomially on the
$c_{ij}^k$'s.  This defines a polynomial function $h:V_{4,2}\longrightarrow\RR$,
which turns to be $\Sl_4(\RR)\times\Sl_2(\RR)$-invariant.

It is not hard to see that $h(\lb)\ne 0$ if and only if $\lb$ is isomorphic to
either $\hg_3\otimes\CC$ ($h>0$) or $\hg_3\oplus\hg_3$ ($h<0$).  This implies that
the union of the two $\Gl_4(\RR)\times\Gl_2(\RR)$-orbits corresponding to
$\hg_3\oplus\CC$ and $\hg_3\oplus\hg_3$, which coincides with the set of all
Einstein nilradicals of eigenvalue type $(1<2;4,2)$ in $V_{4,2}$, is open and dense
in $V_{4,2}$.  However, recall that the net probability of being an Einstein
nilradical of eigenvalue type $(1<2;4,2)$ in $V_{4,2}$ is $\tfrac{2}{7}$.

One may try to avoid this by working on the quotient space
$V_{q,p}/\Gl_q(\RR)\times\Gl_p(\RR)$, where Theorem \ref{generic} is by the way also
true, but the topology here is so ugly that an open and dense subset can never be
taken as a probability one subset.  In fact, there could be a single point set which
is open and dense.  On the other hand, the coset of $0$ is always in the closure of
any other subset, which shows that this quotient space is far from being $T_1$.

It has very recently appeared in \cite{Nkl4} a complete classification for $2$-step
Einstein nilradicals of type $(2,q)$ for any $q$.  In \cite{Jbl}, a construction
called concatenation of $2$-step nilpotent Lie algebras is used to obtain Einstein
nilradicals of type $(1<2;q,p)$ from smaller ones, as well as many new examples of
$2$-step non-Einstein nilradicals.

\section{Known examples and non examples}\label{ene}

As far as we know, the following is a complete chronological list of nilpotent Lie
algebras which are known to be Einstein nilradicals,  or equivalently, of known
examples of rank-one Einstein solvmanifolds:
\begin{itemize}

\item[ ]
\item \cite{Cart} The Lie algebra of an Iwasawa $N$-group: $G/K$ irreducible symmetric space of
noncompact type and $G=KAN$ the Iwasawa decomposition.

\item[ ]
\item\cite{GndPttVnb} Nilradicals of normal $j$-algebras (i.e. of
noncompact homogeneous K$\ddot{{\rm a}}$hler Einstein spaces).

\item[ ]
\item\cite{Alk, Crt} Nilradicals of homogeneous quaternionic K$\ddot{{\rm
a}}$hler spaces.

\item[ ]
\item\cite{Dlf} Certain $2$-step nilpotent Lie algebras for which there
is a basis with very uniform properties (see also \cite[1.9]{Wlt}).

\item[ ]
\item\cite{Bgg} $H$-type Lie algebras (see also \cite{Lnz}).

\item[ ]
\item\cite{EbrHbr, manus} Nilpotent Lie algebras admitting a naturally reductive left invariant metric.

\item[ ]
\item\cite{Hbr} Families of deformations of Lie algebras of Iwasawa $N$-groups in the rank-one
case.

\item[ ]
\item\cite{Fan1,Fan2} Certain $2$-step nilpotent Lie algebras
constructed via Clifford modules.

\item[ ]
\item\cite{GrdKrr} A $2$-parameter family of $2$-step nilpotent Lie algebras of type $(3,6)$ and certain
modifications of the Lie algebras of Iwasawa $N$-groups (rank $\geq 2$).

\item[ ]
\item\cite{finding} Any nilpotent Lie algebra with a codimension one abelian ideal.

\item[ ]
\item\cite{finding} A curve of $6$-step nilpotent Lie algebras of dimension $7$, which
is the lowest possible dimension for a continuous family.

\item[ ]
\item\cite{Mor} (and Yamada), Certain $2$-step nilpotent Lie algebras defined from
subsets of fundamental roots of complex simple Lie algebras.

\item[ ]
\item\cite{finding} Any nilpotent Lie algebra of dimension $\leq 5$.

\item[ ]
\item\cite{Wll} Any nilpotent Lie algebra of dimension $6$.

\item[ ]
\item\cite{inter} A curve of $2$-step nilpotent Lie algebras of type $(5,5)$.

\item[ ]
\item\cite{Krr} A $2$-parameter family of deformations of the nilradical of the $12$-dimensional quaternionic hyperbolic space.

\item[ ]
\item\cite{Pyn} Any filiform (i.e. $n$-dimensional and $(n-1)$-step nilpotent) Lie algebra with at least two linearly independent semisimple derivations.

\item[ ]
\item\cite{einsteinsolv} Certain $2$-step nilpotent Lie algebras attached to graphs as soon as a uniquely defined weighting on the graph is positive.  Regular graphs and trees without any edge adjacent to four or more edges are positive.

\item[ ]
\item\cite{Nkl1} The free $p$-step nilpotent Lie algebras $\fg(m,p)$ on $m$ generators for $p=1,2$; $p=3$ and $m=2,3,4,5$; $p=4$ and $m=2$; $p=5$ and $m=2$.

\item[ ]
\item\cite{Nkl2} Several families of filiform Lie algebras.

\item[ ]
\item\cite{Tmr} The nilradical of any parabolic subalgebra of a semisimple
Lie algebra.

\item[ ]
\item\cite{Nkl3} Any $2$-step nilpotent Lie algebra of type $(p,q)$ (i.e. $p+q$-dimensional and $p$-dimensional derived algebra) with $q\leq 5$ and $(p,q)\ne (5,5)$, with the only exceptions of the real forms of six complex algebras of type $(6,5)$ and three of type $(7,5)$.

\item[ ]
\end{itemize}

We now give an up to date list of $\NN$-graded nilpotent Lie
algebras which are not Einstein nilradicals, that is, they do not admit any
nilsoliton metric.

\begin{itemize}
\item[ ]
\item\cite{einsteinsolv} Three $6$-step nilpotent Lie algebras of dimension $7$, and
certain $2$-step nilpotent Lie algebras attached to graphs in any dimension $\geq
11$ (only finitely many in each dimension).

\item[ ]
\item\cite{Nkl1} The free $p$-step nilpotent Lie algebras $\fg(m,p)$ on $m$ generators
for $p=3$ and $m\geq 6$; $p=4$ and $m\geq 3$; $p=5$ and $m\geq 3$; $p\geq 6$.

\item[ ]
\item\cite{Nkl2} Many filiform Lie algebras starting from dimension $8$ (see also \cite{Arr}).

\item[ ]
\item\cite{Nkl3} Real forms of six complex $2$-step nilpotent Lie algebras of
type $(6,5)$ and three of type $(7,5)$.

\item[ ]
\item\cite{Wll2} Two curves of $2$-step nilpotent Lie algebras of type $(3,6)$.

\item[ ]
\end{itemize}

\section{A stratification for the variety of nilpotent Lie algebras}\label{st}

In this section, we define a $\G$-invariant stratification for the representation
$V=\lam$ of $\G$ by adapting to this context the construction given in \cite[Section
12]{Krw1} for reductive group representations over an algebraically closed field.
This construction, in turn, is based on some instability results proved in
\cite{Kmp} and \cite{Hss}.  We decided to give in \cite[Section 2]{standard} a
self-contained proof of all these results, bearing in mind that a direct application
of them does not seem feasible (see also \cite{strata}).

We shall use the notation given in Section \ref{Tb}.  For any $\mu\in V$ we have
that
$$
\lim\limits_{t\to\infty}e^{tI}.\mu= \lim\limits_{t\to\infty}e^{-t}\mu=0,
$$
and hence $0\in\overline{\G.\mu}$, that is, any element of $V$ is unstable for our
$\G$-action (see Appendix).  Therefore, in order to distinguish two elements of $V$
from the point of view of geometric invariant theory, we would need to measure in
some sense `how' unstable each element of $V$ is.  Maybe the above is not the
optimal way to go to $0$ along the orbit starting from $\mu$.

Let us consider $\mu\in V$ and $\alpha\in\dca$, where $\dca$ denotes the set of all
$n\times n$ matrices which are diagonalizable, that is,
$$
\dca=\bigcup_{g\in\G}g\tg g^{-1}.
$$
Thus $\pi(\alpha)$ is also diagonalizable (see (\ref{actiong})), say with
eigenvalues $a_1,...,a_r$ and eigenspace decomposition  $V=V_1\oplus...\oplus V_r$.
This implies that if $\mu\ne 0$ and $\mu=\mu_1+...+\mu_r$, $\mu_i\in V_i$, then
$$
e^{-t\alpha}.\mu= \sum_{i=1}^r e^{-ta_i}\mu_i,
$$
and so $e^{-t\alpha}.\mu$ goes to $0$ when $t\to\infty$ if and only if $\mu_i=0$ as
soon as $a_i\leq 0$.  Moreover, in that case, the positive number
$$
m(\mu,\alpha):=\min\{ a_i : \mu_i\ne 0\},
$$
measures the degree of instability of $\mu$ relative to $\alpha$, in the sense that
the train has not arrived until the last wagon has.  Indeed, the larger
$m(\mu,\alpha)$ is, the faster $e^{-t\alpha}.\mu$ will converge to $0$ when
$t\to\infty$.  Recall that for an action in general the existence of such $\alpha$
for any unstable element is guaranteed by Theorem \ref{RS}, (iv).

Notice that $m(\mu,c\alpha)=c m(\mu,\alpha)$ for any $c>0$.  We can therefore
consider the most efficient directions (up to the natural normalization) for a given
$\mu\in V$, given by
$$
\Lambda(\mu):=\left\{\beta\in\dca: m(\mu,\beta)=1=\sup\limits_{\alpha\in\dca}\left\{
m(\mu,\alpha):\tr{\alpha^2}=\tr{\beta^2}\right\}\right\}.
$$
A remarkable fact is that $\Lambda(\mu)$ lie in a single conjugacy class, that is,
there exists an essentially unique direction which is `most responsible' for the
instability of $\mu$.  All the parabolic subgroups $P_{\beta}$ of $\G$
naturally associated to any $\beta\in\Lambda(\mu)$ defined in (\ref{para}) coincide, and
hence they define a unique parabolic subgroup $P_{\mu}$ which acts transitively on
$\Lambda(\mu)$ by conjugation. A very nice property $P_{\mu}$ has is that
\begin{equation}\label{aut}
\Aut(\mu)\subset P_{\mu}.
\end{equation}

Since
$$
\Lambda(g.\mu)=g\Lambda(\mu)g^{-1}, \qquad \forall\mu\in V, \quad g\in\G,
$$
we obtain that $\Lambda(g.\mu)$ will meet the Weyl chamber $\tg^+$ for some
$g\in\G$, and the intersection set will consist of a single element $\beta\in\tg^+$
(see (\ref{weyl})).

Summarizing, we have been able to attach to each nonzero $\mu\in V$, and actually to
each nonzero $\G$-orbit in $V$, a uniquely defined $\beta\in\tg^+$ which comes from
instability considerations.

\begin{definition}
Under the above conditions, we say that $\mu\in\sca_{\beta}$ and call the subset
$\sca_{\beta}\subset V$ a {\it stratum}.
\end{definition}

We note that $\sca_{\beta}$ is $\G$-invariant for any $\beta\in\tg^+$ and
$$
V\smallsetminus\{ 0\}=\bigcup\limits_{\beta\in\tg^+}\sca_{\beta},
$$
a disjoint union.  An alternative way to define $\sca_{\beta}$ is
$$
\sca_{\beta}=\G.\left\{\mu\in V:\tfrac{\beta}{||\beta||^2}\in\Lambda(\mu)\right\},
$$
which actually works for any $\beta\in\tg$.  From now on, we will always denote by
$\mu_{ij}^k$ the structure constants of a vector $\mu\in V$ with respect to the
basis $\{ v_{ijk}\}$:
$$
\mu=\sum\mu_{ij}^kv_{ijk}, \qquad \mu_{ij}^k\in\RR, \qquad {\rm i.e.}\quad
\mu(e_i,e_j)=\sum_{k=1}^n\mu_{ij}^ke_k, \quad i<j.
$$
Each nonzero $\mu\in V$ uniquely determines an element $\beta_{\mu}\in\tg$ given by
$$
\beta_{\mu}=\mcc\left\{\alpha_{ij}^k:\mu_{ij}^k\ne 0\right\}.
$$
Recall that $\mcc(X)$ denotes the unique element of minimal norm in the convex hull
$\CH(X)$ of a subset $X\subset\tg$, and thus $\beta_{\mu}$ has rational coefficients
(see (\ref{rat})).  We also note that $\beta_{\mu}$ is always nonzero since
$\tr{\alpha_{ij}^k}=-1$ for all $i<j$ and consequently $\tr{\beta_{\mu}}=-1$. If for
$\mu\in V$ we define $\Lambda_T(\mu)$ as above but by replacing $\dca$ with the set
of diagonal matrices $\tg$, then one can prove that
$$
\Lambda_T(\mu)=\left\{\tfrac{\beta_{\mu}}{||\beta_{\mu}||^2}\right\},
$$
that is, $\beta_{\mu}$ is the (unique) `most responsible' direction for the
instability of $\mu$ with respect to the action of the torus $T$ with Lie algebra
$\tg$ on $V$.  Another equivalent definition for the stratum $\sca_{\beta}$,
$\beta\in\tg$, is given by
$$
\sca_{\beta}=\Big\{\mu\in V\smallsetminus\{ 0\}:\beta\;\mbox{is an element of
maximal norm in}\;\{\beta_{g.\mu}:g\in\G\}\Big\}.
$$
If $\mu$ runs through $V$, there are only finitely many possible vectors
$\beta_{\mu}$, and consequently the set $\{\beta\in\tg:\sca_{\beta}\ne\emptyset\}$
is finite.  We furthermore get from this new description that if $\beta\in\tg$
satisfies $\sca_{\beta}\ne \emptyset$ then $\beta$ has rational coefficients and
\begin{equation}\label{trb}
\tr{\beta}=-1.
\end{equation}

\begin{remark}
A very illustrative exercise is to consider the action given in Example
\ref{hompol1}, draw the nice picture of its weights, detect all possible
$\beta_{p}$'s and try to figure out which of them actually determine a nonempty
stratum (i.e. $\sca_{\beta_p}\ne\emptyset$).
\end{remark}

Recall from Section \ref{va} that the moment map $m$ for the $\G$-representation $V$
plays a fundamental role in the study of Einstein solvmanifolds and nilsolitons, as
$m(\mu)=-\tfrac{1}{\scalar(\mu)}\Ric_{\mu}$, where $\Ric_{\mu}$ and $\scalar(\mu)$
denote the Ricci operator and the scalar curvature of $(N_{\mu},\ip)$, respectively.
The square norm functional $F(\mu)=||m(\mu)||^2$ therefore provides a natural
curvature functional on the space $\nca$ of all left invariant metrics on
$n$-dimensional nilpotent Lie groups whose critical points are precisely nilsoliton
metrics (see Theorem \ref{crit}).

We have collected in the following theorem some relationships between $m$, $F$ and
the strata. Let $p_{\tg}(\alpha)$ denote the orthogonal projection on $\tg$ of an
$\alpha\in\sym(n)$ (i.e. the diagonal part of $\alpha$).

\begin{theorem}\label{strataF}\cite{einsteinsolv}
Let $\mu=\sum\mu_{ij}^kv_{ijk}$ be a nonzero element of $V$.
\begin{itemize}
\item[(i)] $p_{\tg}(m(\mu)) =\tfrac{2}{||\mu||^2}\sum\limits_{i<j}(\mu_{ij}^k)^2
\alpha_{ij}^k\in\CH\left\{\alpha_{ij}^k:\mu_{ij}^k\ne 0\right\}$.

\item[(ii)] $F(\mu)\geq ||\beta||^2$ for any $\mu\in\sca_{\beta}$.

\item[(iii)] If $\inf{F(\G.\mu)}=||\beta_{\mu}||^2$ then $\mu\in\sca_{\beta_{\mu}}$.

\item[(iv)] If $\mu\in\nca$, $m(\mu)\in\tg$ and $S_{\mu}$ is Einstein then $\mu\in\sca_{m(\mu)}$.

\item[(v)] For $\mu\in\sca_{\beta}\cap\nca$, the following conditions are equivalent:
\begin{itemize}
\item[(a)] $S_{\mu}$ is Einstein.

\item[(b)] $S_{\mu}$ is Einstein of eigenvalue-type $\beta+||\beta||^2I$ (up to a positive multiple).

\item[(c)] $m(\mu)$ is conjugate to $\beta$.

\item[(d)] $F(\mu)=||\beta||^2$.
\end{itemize}
\end{itemize}
\end{theorem}

It follows from part (v) in the above theorem that the stratum $\sca_{\beta}$ to
which $\mu$ belongs determines the eigenvalue type of a potential Einstein
solvmanifold $S_{g.\mu}$, $g\in\G$ (if any), and so the stratification provides a
convenient tool to produce existence results as well as obstructions for nilpotent
Lie algebras to be an Einstein nilradical.  Thus $\beta$ plays a role similar to
the one played by the pre-Einstein derivation $\phi$ (see Definition \ref{preE}).  The subtle
relationship between $\beta$ and $\phi$ will be explained in Section \ref{smadre}.

We will now give a description of the strata in terms of semistable vectors (see
Appendix).  For each $\beta\in\tg$ consider the sets
$$
\begin{array}{l}
Z_{\beta}=\{ \mu\in V: \la\beta,\alpha_{ij}^k\ra=||\beta||^2,\quad\forall \mu_{ij}^k\ne 0\}, \\ \\

W_{\beta}=\{\mu\in
V:\la\beta,\alpha_{ij}^k\ra\geq||\beta||^2,\quad\forall\mu_{ij}^k\ne 0\}. \\ \\

Y_{\beta}=\{ \mu\in W_{\beta}:\la\beta,\alpha_{ij}^k\ra=||\beta||^2,\quad\mbox{for
at least one}\; \mu_{ij}^k\ne 0\}.
\end{array}
$$

Notice that $Z_{\beta}$ is actually the eigenspace of $\pi(\beta)$ with eigenvalue
$||\beta||^2$, and so $\mu\in Z_{\beta}$ if and only if
$\beta+||\beta||^2I\in\Der(\mu)$.  We also note that $W_{\beta}$ is the direct sum
of all the eigenspaces of $\pi(\beta)$ with eigenvalues $\geq ||\beta||^2$, and
since $Z_{\beta}\subset Y_{\beta}\subset W_{\beta}$, they are all
$\G_{\beta}$-invariant, where $\G_{\beta}$ is the centralizer of $\beta$ in $\G$.

Let $\g_{\beta}$ denote the Lie algebra of $\G_{\beta}$, that is,
$$
\g_{\beta}=\{\alpha\in\g:[\alpha,\beta]=0\},
$$
and let $G_{\beta}$ be any reductive subgroup of $\G$ with Lie algebra
$\ggo_{\beta}$, the orthogonal complement of $\beta$ in $\g_{\beta}$.
Thus
\begin{equation}\label{descgb}
\g_{\beta}=\ggo_{\beta}\oplus\RR\beta,
\end{equation}

\noindent is an orthogonal decomposition,
$$
\ggo_{\beta}=(\sog(n)\cap\ggo_{\beta})\oplus(\sym(n)\cap\ggo_{\beta})
$$
is a Cartan decomposition and
$\tg\cap\ggo_{\beta}=\{\alpha\in\tg:\la\alpha,\beta\ra=0\}$ is a maximal abelian
subalgebra of $\sym(n)\cap\ggo_{\beta}$.  Recall that if $\sca_{\beta}\ne\emptyset$
then $\beta$ is rational and so such a reductive group $G_{\beta}$ does exist.

\begin{definition}
A vector $\mu\in V$ is called $G_{\beta}$-{\it semistable} if
$0\notin\overline{G_{\beta}.\mu}$.
\end{definition}

\begin{theorem}\label{zybeta}\cite{standard}
For any $\beta\in\tg$, the $\G_{\beta}$-invariant subsets
$Z_{\beta}^{ss}:=Z_{\beta}\cap\sca_{\beta}$ and
$Y_{\beta}^{ss}:=Y_{\beta}\cap\sca_{\beta}$ satisfy:
\begin{itemize}
\item[(i)] $\sca_{\beta}=\Or(n).Y_{\beta}^{ss}$.

\item[(ii)] $Y_{\beta}^{ss}=\{ \mu\in\sca_{\beta}:\beta_{\mu}=\beta\}$.

\item[(iii)] $Z_{\beta}^{ss}$ is the set of $G_{\beta}$-semistable vectors in $Z_{\beta}$.

\item[(iv)] $Y_{\beta}^{ss}$ is the set of $G_{\beta}$-semistable vectors in $W_{\beta}$.
\end{itemize}
\end{theorem}

We summarize in the following theorem the main properties of the $\G$-invariant
stratification of the vector space $V$ given above.

\begin{theorem}\label{strata}\cite{standard}
There exists a finite subset $\bca\subset\tg^+$, and for each $\beta\in\bca$ a
$\G$-invariant subset $\sca_{\beta}\subset V$ (a {\it stratum}) such that
$$
V\smallsetminus\{ 0\}=\bigcup_{\beta\in\bca}\sca_{\beta} \qquad \mbox{(disjoint
union)}.
$$
If $\mu\in\sca_{\beta}$ then
\begin{equation}\label{adbeta}
\left\la[\beta,D],D\right\ra\geq 0 \qquad\forall\; D\in\Der(\mu)
\qquad(\mbox{equality holds}\;\Leftrightarrow [\beta,D]=0)
\end{equation}

\noindent and
\begin{equation}\label{betapos}
\beta+||\beta||^2I \quad\mbox{is positive definite for all}\; \beta\in\bca
\;\mbox{such that}\; \sca_{\beta}\cap\nca\ne\emptyset.
\end{equation}

If in addition $\mu\in Y_{\beta}^{ss}$, i.e.
\begin{equation}\label{cond}
\min\left\{\la\beta,\alpha_{ij}^k\ra:\mu_{ij}^k\ne 0\right\}=||\beta||^2, \quad
\mbox{or equivalently}\quad 0\notin\overline{G_{\beta}.\mu},
\end{equation}

\noindent then
\begin{equation}\label{betaort}
\tr{\beta D}=0 \quad\forall\; D\in\Der(\mu),
\end{equation}

\noindent and
\begin{equation}\label{delta}
\left\la\pi\left(\beta+||\beta||^2I\right)\mu,\mu\right\ra\geq 0 \qquad
(\mbox{equality holds}\;\Leftrightarrow \beta+||\beta||^2I\in\Der(\mu)).
\end{equation}

Moreover, condition {\rm (\ref{cond})} is always satisfied by some $g.\mu$ with
$g\in\Or(n)$, for any $\mu\in\sca_{\beta}$.
\end{theorem}

\begin{remark} We note that (\ref{betapos}) is actually the only result stated in this section where we really need $\mu$ to be
a nilpotent Lie algebra, and not just any vector in $V$.  It is known for instance
that semisimple Lie algebras lie in the stratum $\sca_{\beta}$ for
$\beta=-\tfrac{1}{n}I$, and consequently $\beta+||\beta||^2I=0$ (see \cite{strata}).
\end{remark}

\section{The stratification and the standard condition}\label{proof}

We now apply the stratification described in Section \ref{st} to prove that Einstein
solvmanifolds are all standard.

Let $S$ be a solvmanifold, that is, a simply connected solvable Lie group endowed
with a left invariant Riemannian metric.  Let $\sg$ be the Lie algebra of $S$ and
let $\ip$ denote the inner product on $\sg$ determined by the metric.  We consider
the orthogonal decomposition $\sg=\ag\oplus\ngo$, where $\ngo=[\sg,\sg]$. Recall
that $S$ is called standard if $[\ag,\ag]=0$.  The mean curvature vector of $S$ is
the only element $H\in\ag$ which satisfies $\la H,A\ra=\tr{\ad{A}}$ for any
$A\in\ag$.  If $B$ denotes the symmetric map defined by the Killing form of $\sg$
relative to $\ip$ then $B(\ag)\subset\ag$ and $B|_{\ngo}=0$ as $\ngo$ is contained
in the nilradical of $\sg$.  The Ricci operator $\Ricci$ of $S$ is given by (see for
instance \cite[7.38]{Bss}):
\begin{equation}\label{ricci}
\Ricci=R-\unm B-S(\ad{H}),
\end{equation}

\noindent where $S(\ad{H})=\unm(\ad{H}+(\ad{H})^t)$ is the symmetric part of
$\ad{H}$ and $R$ is the symmetric operator defined by
\begin{equation}\label{R}
\la Rx,y\ra=-\unm\displaystyle{\sum\limits_{ij}}\la [x,x_i],x_j\ra\la [y,x_i],x_j\ra
+\unc\displaystyle{\sum\limits_{ij}}\la [x_i,x_j],x\ra\la [x_i,x_j],y\ra,
\end{equation}

\noindent for all $x,y\in\sg$, where $\{ x_i\}$ is any orthonormal basis of
$(\sg,\ip)$.

It is proved in \cite[Propositions 3.5, 4.2]{minimal} that $R$ is the only symmetric
operator on $\sg$ such that
\begin{equation}\label{Rmm}
\tr{RE}=\unc\la \pi(E)[\cdot,\cdot],[\cdot,\cdot]\ra, \qquad\forall E\in\End(\sg),
\end{equation}

\noindent where we are considering $[\cdot,\cdot]$ as a vector in
$\Lambda^2\sg^*\otimes\sg$, $\ip$ is the inner product defined in (\ref{innV}) and
$\pi$ is the representation given in (\ref{actiong}) (see the notation in Section
\ref{Tb} and replace $\RR^n$ by $\sg$).  This is equivalent to saying that
$$
m([\cdot,\cdot])=\tfrac{4}{|| [\cdot,\cdot]||^2}R,
$$
where $m:\Lambda^2\sg^*\otimes\sg\longrightarrow\sym(\sg)$ is the moment map for the
action of $\Gl(\sg)$ on $\Lambda^2\sg^*\otimes\sg$ (see (\ref{defmm})).  Thus the
`anonymous' tensor $R$ in formula (\ref{ricci}) for the Ricci operator is precisely
the value of the moment map at the Lie bracket $[\cdot,\cdot]$ of $\sg$ (up to
scaling).

We therefore obtain from (\ref{ricci}) and (\ref{Rmm}) that $S$ is an Einstein
solvmanifold with $\Ricci=cI$, if and only if, for any $E\in\End(\sg)$,
\begin{equation}\label{einstein}
\tr{\left(cI+\unm B+S(\ad{H})\right)E}= \unc\la
\pi(E)[\cdot,\cdot],[\cdot,\cdot]\ra.
\end{equation}

Let $S$ be an Einstein solvmanifold with $\Ricci=cI$.  We can assume that $S$ is not
unimodular by using \cite{Dtt}, thus $H\ne 0$ and $\tr{\ad{H}}=||H||^2>0$.  By
letting $E=\ad{H}$ in (\ref{einstein}) we get
\begin{equation}\label{c}
c=-\tfrac{\tr{S(\ad{H})^2}}{\tr{S(\ad{H})}}<0.
\end{equation}

In order to apply the results in Section \ref{st}, we identify $\ngo$ with $\RR^n$
via an orthonormal basis $\{ e_1,...,e_n\}$ of $\ngo$ and we set
$\mu:=[\cdot,\cdot]|_{\ngo\times\ngo}$.  In this way, $\mu$ can be viewed as an
element of $\nca\subset V$. If $\mu\ne 0$ then $\mu$ lies in a unique stratum
$\sca_{\beta}$, $\beta\in\bca$, by Theorem \ref{strata}, and it is easy to see that
we can assume (up to isometry) that $\mu$ satisfies (\ref{cond}), so that one can
use all the additional properties stated in the theorem.  In particular, the
following crucial technical result follows.  Consider $E_{\beta}\in\End(\sg)$
defined by
$$
E_{\beta}=\left[\begin{smallmatrix} 0&0\\
0&\beta+||\beta||^2I\end{smallmatrix}\right],
$$
that is, $E|_{\ag}=0$ and $E|_{\ngo}=\beta+||\beta||^2I$.

\begin{lemma}\label{pie}
If $\mu\in\sca_{\beta}$ satisfies {\rm (\ref{cond})} then
$\la\pi(E_{\beta})[\cdot,\cdot],[\cdot,\cdot]\ra\geq 0$.
\end{lemma}

We then apply (\ref{einstein}) to $E_{\beta}\in\End(\sg)$ and obtain from Lemma
\ref{pie}, (\ref{c}), (\ref{trb}) and (\ref{betaort}) that
$$
\tr{S(\ad{H})^2}\tr{E_{\beta}^2}\leq (\tr{S(\ad{H})E_{\beta}})^2,
$$
a `backwards' Cauchy-Schwartz inequality.  This turns all inequalities which
appeared in the proof of Lemma \ref{pie} into equalities, in particular:
$$
\unc\sum_{rs}\la(\beta+||\beta||^2I)[A_r,A_s],[A_r,A_s]\ra=0,
$$
where $\{ A_i\}$ is an orthonormal basis of $\ag$.  We finally get that $\ag$ is
abelian since $\beta+||\beta||^2I$ is positive definite by (\ref{betapos}).

\section{The stratification and Einstein solvmanifolds via closed orbits}\label{smadre}

We shall describe in this section some other applications of the strata defined in
Section \ref{st} to the study of Einstein solvmanifolds.

Let $\ngo$ be a nonabelian nilpotent Lie algebra of dimension $n$.  We fix any basis
$\{ X_1,...,X_n\}$ of $\ngo$ and consider the corresponding structural constants:
$$
[X_i,X_j]=\sum_{k=1}^n c_{ij}^kX_k, \qquad 1\leq i<j\leq n.
$$
Let $\beta$ denote the unique element of minimal norm in the convex hull of the set
$\{\alpha_{ij}^k:c_{ij}^k\ne 0\}$, where $\alpha_{ij}^k$ is the diagonal $n\times n$
matrix $-E_{ii}-E_{jj}+E_{kk}$.  Notice that $\tr{\beta}=-1$, and so $\beta$ is
always nonzero.  We define the Lie algebra
$$
\ggo_{\beta}=\{\alpha\in\End(\ngo):[\alpha,\beta]=0, \quad\tr{\alpha\beta}=0\},
$$
and take any reductive subgroup $G_{\beta}$ of $\Gl(\ngo)$ with Lie algebra
$\ggo_{\beta}$ (existence is guaranteed by rationality of $\beta$).

Recall that the Lie bracket $\lb$ of $\ngo$ belongs to the vector space
$\Lambda^2\ngo^*\otimes\ngo$ of skew-symmetric bilinear maps from $\ngo\times\ngo$
to $\ngo$, on which $\Gl(\ngo)$ is acting naturally by
$g.\lb=g[g^{-1}\cdot,g^{-1}\cdot]$.

\begin{theorem}\label{madre}\cite{einsteinclosed}
Let $\ngo$ be a nonabelian nilpotent Lie algebra and for any basis $\{
X_1,...,X_n\}$ of $\ngo$ consider $\beta$ and $G_{\beta}\subset\Gl(\ngo)$ as defined
above.

\begin{itemize}
\item[(i)] If the orbit $G_{\beta}.\lb$ is closed in $\Lambda^2\ngo^*\otimes\ngo$ then $\ngo$ is an Einstein nilradical and $\beta+||\beta||^2I\in\Der(\ngo)$.

\item[(ii)] If $\ngo$ is an Einstein nilradical, $\beta+||\beta||^2I\in\Der(\ngo)$ and $0\notin\overline{G_{\beta}.\lb}$, then the orbit $G_{\beta}.\lb$ is closed in $\Lambda^2\ngo^*\otimes\ngo$ .
\end{itemize}
\end{theorem}

\begin{remark}
The following example shows that condition $\beta+||\beta||^2I\in\Der(\ngo)$ is
necessary in part (ii) of Theorem \ref{madre}.  Let $\ngo$ be the $4$-dimensional
$3$-step nilpotent Lie algebra with Lie bracket given by
$$
[X_1,X_2]=X_3+X_4, \qquad [X_1,X_3]=X_4.
$$
It is easy to see that $\beta=(-1,-\unm,0,\unm)$ and
$0\notin\overline{G_{\beta}.\lb}$.  If $\lambda$ is defined by
$$
\lambda(X_1,X_2)=X_3, \qquad \lambda(X_1,X_3)=X_4,
$$
then $\lambda\in\Gl_4(\RR).\lb$, $m(\lambda)=\beta$ and
$\beta+||\beta||^2I\in\Der(\lambda)$, from which follows that $\lambda$ is a
nilsoliton and so $\ngo$ is an Einstein nilradical.  However,
$$
\lambda=\lim_{t\to\infty} e^{-t\alpha}.\lb\in\overline{G_{\beta}.\lb}, \qquad
\mbox{for}\quad\alpha=(1,0,1,2),
$$
and thus $G_{\beta}.\lb$ is not closed.  Indeed, $\lambda\notin G_{\beta}.\lb$ since
$\beta+||\beta||^2I\in\Der(\lambda)$ and $\beta+||\beta||^2I\notin\Der(\ngo)$.
\end{remark}

Recall that $\beta$ has entries in $\QQ$ and so if $\beta\in\tg^+$ and has
eigenvalues $b_1<...<b_r$ with multiplicities $n_1,...,n_r$, respectively, then one can for instance take the reductive group $G_{\beta}$ given by
$$
G_{\beta}=\left\{\left[\begin{smallmatrix} g_1&& \\ &\ddots&\\
&&g_r\end{smallmatrix}\right]:\det{g_1}^{mb_1}...\det{g_r}^{mb_r}=1, \quad
g_i\in\Gl_{n_i}(\RR)\right\},
$$
where $m$ is the least common multiple of the denominators of the $b_i$'s.  If
$\sigma$ is a permutation of $\{ 1,...,n\}$ then the new basis $\{
X_{\sigma(1)},...,X_{\sigma(n)}\}$ of $\ngo$ has structural constants
$$
[X_{\sigma(i)},X_{\sigma(j)}]=\sum_{k=1}^n
c_{\sigma(i)\sigma(j)}^{\sigma(k)}X_{\sigma(k)}, \qquad 1\leq i<j\leq n,
$$
and so the new $\beta$ has eigenvalues $b_{\sigma^{-1}(1)},...,b_{\sigma^{-1}(n)}$
with respective eigenvectors $X_{\sigma(1)},...,X_{\sigma(n)}$ (see the beginning of
the proof of \cite[Theorem 2.10]{standard}).  Therefore, we can always assume that
$\beta\in\tg^+$, up to just a permutation of the basis $\{ X_i\}$.  Otherwise, if
one insists on keeping the original basis, one may take as $G_{\beta}$ the group
$h^{-1}G_{h\beta h^{-1}}h$, where $h\in\Gl(\ngo)$ is a permutation matrix such that
$h\beta h^{-1}\in\tg^+$.

The following two results show the interplay between the stratum $\beta$ and the pre-Einstein derivation $\phi$ (see Definition \ref{preE}), providing in particular a new method to compute $\phi$.

\begin{lemma}\label{madre2}\cite{einsteinclosed}
Let $\ngo$ be a nonabelian nilpotent Lie algebra and for any basis $\{
X_1,...,X_n\}$ of $\ngo$ consider $\beta$ and $G_{\beta}\subset\Gl(\ngo)$ as defined
above.

\begin{itemize}
\item[(i)] $\lb\in\sca_{\beta}$ if and only if $0\notin\overline{G_{\beta}.\lb}$.

\item[(ii)] If $\beta+||\beta||^2I\in\Der(\ngo)$ and $0\notin\overline{G_{\beta}.\lb}$ then
$$
\phi:=\tfrac{1}{||\beta||^2}(\beta+||\beta||^2I)
$$
is a pre-Einstein derivation of $\ngo$ and $G_{\phi}.\lb$ is closed if and only if
$G_{\beta}.\lb$ is closed.
\end{itemize}
\end{lemma}

\begin{remark}
We conclude from Lemma \ref{madre2} that Theorem \ref{madreN} and Theorem
\ref{madre} are equivalent.
\end{remark}

\begin{remark}
The stratum a given nilpotent Lie algebra belongs to provides useful information on
its automorphism group.  Indeed, let $\ngo$ be a nilpotent Lie algebra and for any
basis $\{ X_1,...,X_n\}$ of $\ngo$ consider $\beta$ as defined above. If
$\lb\in\sca_{\beta}$, then $\Aut(\ngo)\subset P_{\beta}$ by (\ref{aut}).
\end{remark}

\begin{lemma}\label{madre3}\cite{einsteinclosed}
Let $\ngo$ be a nonabelian nilpotent Lie algebra and let $\phi$ be a pre-Einstein
derivation of $\ngo$ with basis of eigenvectors $\{ X_1,...,X_n\}$ and define
$$
\beta:=\tfrac{1}{n-\tr{\phi}}(\phi-I).
$$
Then $\lb\in\sca_{\beta}$ if and only if $0\notin\overline{G_{\beta}.\lb}$, and in
that case, $\beta=\mcc(\{\alpha_{ij}^k:c_{ij}^k\ne 0\})$.
\end{lemma}

It follows from Lemma \ref{madre3} that if $\phi$ is a pre-Einstein derivation of
$\ngo$ then $\phi>0$ (see (\ref{betapos})) and $\ad{\phi}\geq 0$ (see
(\ref{adbeta})) are necessary conditions in order to have
$0\notin\overline{G_{\beta}.\lb}$ (i.e. $\lb\in\sca_{\beta}$).  These conditions are not however sufficient (compare with the paragraph below Theorem \ref{eupreE}).  For instance, any free nilpotent Lie algebra which is not an Einstein nilradical provides a counterexample (see \cite[Remark 2]{Nkl3}).

\section{Open problems}\label{op}

Let $\ngo$ be an $\NN$-graded nilpotent Lie algebra.

\begin{enumerate}
\item {\bf Obstructions}. To find algebraic necessary conditions on $\ngo$ to be an Einstein nilradical.

\item {\bf Existence}. Are there algebraic conditions on $\ngo$ which are sufficient to be an Einstein nilradical?

\item Does the assertion `$\ngo$ is an Einstein nilradical' have probability $1$ in some sense?

\item Does the assertion `$\ngo$ is not an Einstein nilradical' have probability $1$ in some sense?

\item Assume that $\ngo$ is an Einstein nilradical with Lie bracket $\mu_0\in\nca$, and consider the flow $\mu(t)$ defined in (\ref{flow}) with $\mu(0)=\mu_0$.  Does $\lambda=\lim\limits_{t\to\infty}\mu(t)$ necessarily belong to $\G.\mu_0$? (this would provide a nice obstruction).

\item To exhibit an explicit example or prove the existence of a nilpotent Lie algebra which does not admit a nice basis (see Definition \ref{nice}).

\item Are there only finitely many $\NN$-graded filiform Lie algebras which are not Einstein nilradicals in each dimension?
\end{enumerate}

\section{Appendix: Real geometric invariant theory}\label{git}

Let $G$ be a real reductive group acting linearly on a finite dimensional real
vector space $V$ via $(g,v)\mapsto g.v$, $g\in G,v\in V$. The precise definition of
our setting is the one considered in \cite{RchSld}. We also refer to \cite{EbrJbl},
where many results from geometric invariant theory are adapted and proved over
$\RR$.

The Lie algebra $\ggo$ of $G$ also acts linearly on $V$ by the derivative of the
above action, which will be
 denoted by $(\alpha,v)\mapsto\pi(\alpha)v$, $\alpha\in\ggo$, $v\in V$. We consider
 a Cartan decomposition
$\ggo=\kg\oplus\pg$, where $\kg$ is the Lie algebra of a maximal compact subgroup
$K$ of $G$. Endow $V$ with a fixed from now on $K$-invariant inner product $\ip$
such that $\pg$ acts by symmetric operators, and endow $\pg$ with an
$\Ad(K)$-invariant inner product $\ipp$.

The function $m:V\smallsetminus\{ 0\}\longrightarrow\pg$ implicitly defined by
$$
(m(v),\alpha)=\isn\la\pi(\alpha)v,v\ra, \qquad \forall\alpha\in\pg,\; v\in V,
$$
is called the {\it moment map} for the representation $V$ of $G$.  Since
$m(cv)=m(v)$ for any nonzero $c\in\RR$, we also may consider the moment map on the
projective space of $V$, $m:\PP V\mapsto\pg$, with the same notation and definition
as above for $m([v])$, $[v]$ the class of $v$ in $\PP V$.  It is easy to see that
$m$ is $K$-equivariant: $m(k.v)=\Ad(k)m(v)$ for all $k\in K$.

In the complex case (i.e. for a complex representation of a complex reductive
algebraic group), under the natural identifications $\pg=\pg^*=(\im\kg)^*=\kg^*$,
the function $m$ is precisely the moment map from symplectic geometry, corresponding
to the Hamiltonian action of $K$ on the symplectic manifold $\PP V$ (see for
instance the survey \cite{Krw2} or \cite[Chapter 8]{Mmf} for further information).
For real actions, this nice interplay with symplectic geometry is lost ($\PP V$
could even be odd dimensional), but the moment map is nevertheless a very natural
object attached to a real representation encoding a lot of information on the
geometry of $G$-orbits and the orbit space $V/G$.

Let $\mca=\mca(G,V)$ denote the set of {\it minimal vectors}, that is,
$$
\mca=\{ v\in V: ||v||\leq ||g.v||\quad\forall g\in G\}.
$$
For each $v\in V$ define
$$
\rho_v:G\mapsto\RR, \qquad \rho_v(g)=||g.v||^2.
$$
In \cite{RchSld}, it is shown that the nice interplay between closed orbits and
minimal vectors discovered in \cite{KmpNss} for actions of complex reductive
algebraic groups, is still valid in the real situation.

\begin{theorem}\cite{RchSld}\label{RS}
Let $V$ be a real representation of a real reductive group $G$, and let $v\in V$.
\begin{itemize}
\item[(i)] The orbit $G.v$ is closed if and only if $G.v$ meets $\mca$.

\item[(ii)] $v\in\mca$ if and only if $\rho_v$ has a critical point at $e\in G$.

\item[(iii)] If $v\in\mca$ then $G.v\cap\mca=K.v$.

\item[(iv)] The closure $\overline{G.v}$ of any orbit $G.v$ always meets $\mca$.  Moreover, there always exists
$\alpha\in\pg$ such that $\lim\limits_{t\to\infty}\exp(-t\alpha).v=w$ exists and
$G.w$ is closed.

\item[(v)] $\overline{G.v}\cap\mca$ is a single $K$-orbit, or in other words, $\overline{G.v}$ contains a unique
closed $G$-orbit.
\end{itemize}
\end{theorem}

As usual in the real case, classical topology of $V$ is always considered rather
than Zarisky topology, unless explicitly indicated.

Let $(\dif\rho_v)_e:\ggo\mapsto\RR$ denote the differential of $\rho_v$ at the
identity $e$ of $G$. It follows from the $K$-invariance of $\ip$ that
$(\dif\rho_v)_e$ vanishes on $\kg$, and so we can assume that
$(\dif\rho_v)_e\in\pg^*$, the vector space of real-valued functionals on $\pg$.  If
we identify $\pg$ and $\pg^*$ by using $\ipp$, then it is easy to see that
$$
m(v)=\tfrac{1}{2||v||^2}(\dif\rho_v)_e.
$$
The moment map at $v$ is therefore an indicator of the behavior of the norm along
the orbit $G.v$ in a neighborhood of $v$.  It follows from Theorem \ref{RS}, (ii)
that
$$
\mca\smallsetminus\{ 0\} = \{ v\in V\smallsetminus\{ 0\}: m(v)=0\}.
$$
Thus if we consider the functional square norm of the moment map
\begin{equation}\label{norm}
F:V\smallsetminus\{ 0\}\mapsto\RR, \qquad  F(v)=||m(v)||^2,
\end{equation}
which is a $4$-degree homogeneous polynomial times $||v||^{-4}$,  $\mca\setminus\{ 0\}$
coincides with the set of zeros of $F$. It then follows from Theorem \ref{RS}, parts
(i) and (iii), that a nonzero orbit $G.v$ is closed if and only if $F(w)=0$ for some
$w\in G.v$, and in that case, the set of zeros of $F|_{G.v}$ coincides with $K.v$.
Recall that $F$ is scaling invariant and so it is actually a function on any sphere
of $V$ or on $\PP V$.

A natural question arises: what is the role played by the remaining critical points
of $F$ (i.e. those for which $F(v)>0$) in the study of the $G$-orbit space of the
action of $G$ on $V$?. This was independently studied in \cite{Krw1} and \cite{Nss}
in the complex case, who have shown that non-minimal critical points still enjoy
most of the nice properties of minimal vectors stated in Theorem \ref{RS}. In the
real case, the analogues of some of these results have been proved in \cite{Mrn}.

We endow $\PP V$ with the Fubini-Study metric defined by $\ip$ and denote by
$x\mapsto \alpha_x$ the vector field on $\PP V$ defined by $\alpha\in\ggo$ via the
action of $G$ on $\PP V$, that is, $\alpha_x=\ddt|_0\exp(t\alpha).x$. We will also
denote by $F$ the functional $F:\PP V\longrightarrow\RR$, $F([v])=||m([v])||^2$.

\begin{lemma}\cite{Mrn}\label{marian2}
The gradient of the functional $F:V\setminus\{ 0\}\longrightarrow\RR$ is given by
$$
\grad(F)_{v}=\tfrac{4}{||v||^2}\Big(\pi(m(v))v-||m(v)||^2v\Big), \qquad v\in
V\setminus\{ 0\},
$$
and for $F:\PP V\longrightarrow\RR$ we have that
$$
\grad(F)_{[v]}=4m([v])_{[v]}, \qquad [v]\in\PP V.
$$
\end{lemma}

Therefore, $v$ is a critical point of $F$ (or equivalently, of $F|_{G.v}$) if and
only if $v$ is an eigenvector of $\pi(m(v))$, and $[v]$ is a critical point of $F$
(or equivalently, of $F|_{G.[v]}$) if and only if $\exp{tm([v])}$ fixes $[v]$.

\begin{theorem}\cite{Mrn}\label{marian}
Let $V$ be a real representation of a real semisimple Lie group $G$.
\begin{itemize}
\item[(i)] If $x\in\PP V$ is a critical point of $F$ then the functional $F|_{G.x}$ attains its minimum value at
$x$.

\item[(ii)] If nonempty, the critical set of $F|_{G.x}$ consists of a single $K$-orbit.
\end{itemize}
\end{theorem}

\begin{definition}\label{stable}
A nonzero vector $v\in V$ is called {\it unstable} if $0\in\overline{G.v}$, and {\it
semistable} otherwise.  If a semistable vector has in addition compact isotropy
subgroup then it is called {\it stable}.
\end{definition}

If the orbit of a nonzero $v\in V$ is closed then $v$ is clearly semistable. More
generally, $v\in V$ is semistable if and only if the unique (up to $K$-action) zero
of $F$ which belongs to $\overline{G.v}$ is a nonzero vector. On the contrary, any
critical point of $F$ which is not a zero of $F$ is unstable.  Indeed, if
$\pi(m(v))v=cv$, $c=||m(v)||^2>0$ (see Lemma \ref{marian2}), then
$$
\lim_{t\to\infty}\exp(-tm(v)).v=\lim_{t\to\infty}e^{-tc}v=0,
$$
and so $0\in\overline{G.v}$.  Thus the study of critical points of $F$ other than
zeroes gives useful information on the orbit space structure of the subset of all
unstable vectors, often called the {\it nullcone} of $V$.

\begin{example}\label{hompol1}
Let us consider the example of $G=\Sl_3(\RR)$ and $V=P_{3,3}(\RR)$, the vector space
of all homogeneous polynomials of degree $3$ on $3$ variables.  The action is given
by a linear change of variables on the left
$$
(g.p)(x_1,x_2,x_3)=p\left(g^{-1}\left[\begin{smallmatrix} x_1\\ x_2\\ x_3
\end{smallmatrix}\right]\right), \qquad \forall g\in\Sl_3(\RR),\quad p\in P_{3,3}(\RR).
$$
It follows that $\ggo=\slg_3(\RR)$, $K=\SO(3)$, $\kg=\sog(3)$ and $\pg=\sym_0(3)$ is
the space of traceless symmetric $3\times 3$ matrices.  As an $\Ad(K)$-invariant
inner product on $\pg$ we take $(\alpha,\beta)=\tr{\alpha\beta}$, and it is easy to
see that the inner product $\ip$ on $V$ for which the basis of monomials
$$
\{ x^D:=x_1^{d_1}x_2^{d_2}x_3^{d_3}:d_1+d_2+d_3=3,\; D=(d_1,d_1,d_3)\}
$$
is orthogonal and
$$
||x^D||^2=d_1!d_2!d_3!, \qquad\forall D=(d_1,d_2,d_3),
$$
satisfies the required conditions.  Let $E_{ij}$ denote as usual the $n\times n$
matrix whose only nonzero coefficient is a $1$ in the entries $ij$.  Since
$$
\pi(E_{ij})p=\ddt|_0p(e^{-tE_{ij}}\cdot )=-x_j\tfrac{\partial p}{\partial x_i},
$$
we obtain that the moment map $m:P_{3,3}(\RR)\longrightarrow\sym_0(3)$ is given by
$$
m(p)=I-\tfrac{1}{||p||^2}\left[\la x_j\tfrac{\partial p}{\partial x_i},p\ra\right].
$$
We are using here that $\la x_j\tfrac{\partial p}{\partial x_i},p\ra=\la
x_i\tfrac{\partial p}{\partial x_j},p\ra$ for all $i,j$.

It is also easy to see that the action of a diagonal matrix $\alpha\in\slg_3(\RR)$
with entries $a_1,a_2,a_3$ is given by
\begin{equation}\label{diagact}
\pi(\alpha)x^D=-\left(\sum_{i=1}^3a_id_i\right)x^D, \qquad\forall D=(d_1,d_2,d_3).
\end{equation}
A first general observation is that any monomial is a critical point of $F$. Indeed,
$$
m(x^D)=\left[\begin{smallmatrix} 1-d_1&&\\ &1-d_2&\\ &&1-d_3
\end{smallmatrix}\right],
$$
and so $x^D$ is an eigenvector of $m(x^D)$ with eigenvalue $F(x^D)=\sum d_i^2-1$
(see Lemma \ref{marian2}).  It follows that $m(p)=0$ for $p=x_1x_2x_3$, that is, $p$
is a minimal vector and its $\Sl_3(\RR)$-orbit is therefore closed.  We also have in
such case that $p_1=p+x_1^3$ is a semistable vector whose orbit is not closed.
Indeed, by acting by diagonal elements with entries $t,\tfrac{1}{t},1$ we get that
$p+t^dx_1^3\in\Sl_3(\RR).p$ for all $t\ne 0$ and so $p\in\overline{\Sl_3(\RR).p_1}$
(recall that $p$ and $p_1$ can never lie in the same orbit since they have
non-isomorphic isotropy subgroups).

For the vector $q=x_1^2x_3+x_1x_2^2$ we have that
$$
m(q)=\left[\begin{smallmatrix} -\tfrac{1}{2}&&\\ &0&\\
&&\tfrac{1}{2}
\end{smallmatrix}\right].
$$
It follows from (\ref{diagact}) that $\pi(m(q))q=\tfrac{1}{2}q$ proving that $q$ is
a critical point of $F$ with critical value $F(q)=\tfrac{1}{2}>0$.  On the other
hand, the family $p_{a,b}=ax_1^2x_3+bx_2^3$, $a,b\ne 0$, lie in a single orbit and
$$
m(p_{a,b})=I-\tfrac{1}{2a^2+6b^2} \left[\begin{smallmatrix} 4a^2&&\\ &18b^2&\\
&&2a^2
\end{smallmatrix}\right].
$$
It is then easy to see by using (\ref{diagact}) that $p_{a,b}$ is a critical point
if and only if $5a^2=27b^2$ and the critical value equals $\tfrac{155}{49}-3$, a
number smaller than $\unm$. In particular, $p_{a,b}$ can not be in the closure of
the orbit of $q$ by Theorem \ref{marian}, (i).
\end{example}

\bibliographystyle{amsalpha}

\end{document}